\documentclass[reqno]{amsproc}
\usepackage{amsmath,amssymb}
\usepackage[T1]{fontenc}
\usepackage{lmodern}
\usepackage{thmtools}
\usepackage{mathtools}
\usepackage[numbers]{natbib}

\usepackage{hyperref}
\usepackage{xurl}
\usepackage{cleveref}

\usepackage[dvipsnames]{xcolor}

\hypersetup{breaklinks=true,
    linkcolor=Bittersweet,
    urlcolor=blue,
    citecolor=blue,
    colorlinks=true,
pdfborder={0 0 0}}

\numberwithin{equation}{section}

\makeatletter
\renewcommand{\bibsection}{\@bibtitlestyle}
\makeatother

\def\R{\mathbb{R}}
\def\E{\mathbb{E}}

\def\eps{\varepsilon}
\DeclareMathOperator{\tr}{tr}

\DeclareMathOperator{\sign}{sign}
\DeclareMathOperator{\rank}{rank}
\DeclareMathOperator{\Rem}{Rem}

\DeclareMathOperator{\loo}{LOO}
\DeclareMathOperator{\alo}{ALO}
\DeclareMathOperator{\diag}{diag}
\DeclareMathOperator{\ran}{ran}
\DeclareMathOperator{\spn}{span}
\DeclareMathOperator*{\argmin}{arg\,min}

\newcommand{\norm}[1]{\lVert#1\rVert}
\newcommand{\lean}[1]{\nolinkurl{#1}}

\declaretheorem[name=Theorem,numberwithin=section]{theorem}
\declaretheorem[name=Proposition,sibling=theorem]{proposition}

\declaretheorem[name=Corollary,sibling=theorem]{corollary}

\declaretheorem[name=Definition,style=definition]{definition}
\declaretheorem[name=Remark,style=remark,numberwithin=section]{remark}
\declaretheorem[name=Example,style=remark,numberwithin=section]{example}

\crefname{assumption}{assumption}{assumptions}

\begin{document}
\title[ALO-CV and mean-field inference II]{Simultaneous analysis of approximate leave-one-out cross-validation
and mean-field inference II: \\
General convex penalties}
\author{Pierre C. Bellec, Rutgers University}

\begin{abstract}
For smooth penalties and group Lasso,
Part~I~\cite{bellec2025simultaneous} analyzed approximate leave-one-out
(ALO) cross-validation by comparing its weights with the
trace correction used in mean-field inference.
For a $C^2$ penalty $R$, the matrix entering these corrections is
$A=[\nabla^2R(\widehat b)+X^TDX]^{-1}$, where $\widehat b$ is the
full-data estimator and $D$ is the diagonal matrix of fitted loss
curvatures; for group Lasso, the penalty and loss curvatures are
first restricted to the active support before taking the
pseudoinverse.

The obstacle to extending this argument to general convex penalties
was the lack of an explicit formula for the matrix $A$ entering
these corrections.
For every finite-valued strongly convex penalty, we construct such
a matrix from ordinary derivatives of the penalty response at the
observed subgradient. Its projection--pseudoinverse formula separates
the response space and penalty curvature, fixed by conditioning,
from the explicit design dependence. This supplies the missing
ingredient to extend the probabilistic argument of Part~I.
The construction recovers the inverse Hessian for smooth penalties
and the active-group inverse for group Lasso, including inactive
boundary groups.

In the proportional regime, for Gaussian designs with general
positive-definite covariance and covariance-scaled strong convexity,
we obtain unconditional bounds comparing both ALO and mean-field
predictions with exact leave-one-out predictions in regression and
single-index models. The expected average squared prediction error is
$O(n^{-1})$, and the expected absolute criterion error is
$O(n^{-1/2})$, under bounded, 1-Lipschitz loss scores and suitable
test functions. The matrix and ALO weights use only observed data;
the mean-field correction additionally uses the design covariance.
The construction, recovery identities and statistical guarantees
are formalized in Lean~4.
\end{abstract}
\maketitle

\tableofcontents

\section{Approximate leave-one-out and mean-field inference}
\label{sec:introduction}

The preceding paper, called \emph{Part~I}~\cite{bellec2025simultaneous} here, established a
probabilistic argument for approximate leave-one-out cross-validation
(ALO) and mean-field inference. Its ALO analysis covered twice
continuously differentiable penalties and group Lasso. The purpose of
Part~II is to remove this restriction: the same conclusions hold for a
general finite-valued strongly convex penalty, without a smoothness
assumption on the penalty.

For the statistical bounds, let $n\ge1$ and assume that the rows of
$X\in\R^{n\times p}$ are independent $N(0,\Sigma)$ vectors with
$\Sigma\succ0$, and that $p\le n/\delta$ for some $\delta>0$.
Let $R$, including its tuning parameters, be deterministic; it may
depend on $n,p$. Assume it is finite valued, minimized at zero, and
strongly convex
in the covariance metric: for some $\mu>0$, the function
\begin{equation}
 b\longmapsto R(b)-\frac{n\mu}{2}\norm{\Sigma^{1/2}b}^2
 \quad\text{is convex}.
 \label{eq:intro-strong-convexity}
\end{equation}
The constants $\mu,\delta$ are kept fixed as $n,p$ vary in the
rate statements below.
Assume the losses $L_y$ are convex and differentiable, with jointly
measurable scores $L_y'$ that are bounded in absolute value by one
and 1-Lipschitz. Define the full-data estimator by
\begin{equation}
 \widehat b=\argmin_{b\in\R^p}
   \left\{\sum_{i=1}^n L_{y_i}(x_i^Tb)+R(b)\right\},
 \label{eq:intro-fit}
\end{equation}
where $x_i^T$ is row $i$ of $X$. Write the fitted predictions,
loss scores and penalty subgradient as
\begin{align}
 t_i&=x_i^T\widehat b,\qquad s_i=L_{y_i}'(t_i),\notag\\
 q&=-X^Ts\in\partial R(\widehat b).
 \label{eq:intro-observed-inputs}
\end{align}
Let $d_i=(L_{y_i}')'(t_i)$ and $D=\diag(d_i)$.
We take $d_i=0$ where this derivative does not exist; $0\le d_i\le1$.

The leave-one-out fit removes observation $i$ while keeping the same
penalty:
\begin{equation}
 \widehat b^{(-i)}
 =\argmin_{b\in\R^p}
   \left\{\sum_{\substack{k=1\\k\ne i}}^nL_{y_k}(x_k^Tb)+R(b)\right\}.
 \label{eq:intro-loo-fit}
\end{equation}
Writing $t_i^{\rm LOO}=x_i^T\widehat b^{(-i)}$, the leave-one-out
criterion for a test function $g$ is
\begin{equation}
 \loo=\frac1n\sum_{i=1}^n g(t_i^{\rm LOO},y_i).
 \label{eq:intro-loo-criterion}
\end{equation}
Computing it requires fitting $n$ estimators. ALO instead approximates
these fits from the single full-data fit $\widehat b$
\cite{rad2018scalable,wang2018approximate_fast_tuning}.

To see where the matrix $A$ comes from, suppose for the moment that
$R$ and the losses are twice continuously differentiable. The
full-data first-order condition is
$X^Ts+\nabla R(\widehat b)=0$, so the gradient of the leave-one-out
objective at $\widehat b$ is $-s_ix_i$. Its Hessian there is
\[
 H_{-i}=\nabla^2R(\widehat b)+\sum_{\substack{k=1\\k\ne i}}^nd_kx_kx_k^T.
\]
One Newton step from $\widehat b$ therefore gives the approximation
\begin{equation}
 \widetilde b^{(-i)}
 =\widehat b+H_{-i}^{-1}x_is_i.
 \label{eq:intro-newton}
\end{equation}
Strong convexity of $R$ makes these Hessians positive definite.
Define the inverse Hessian of the full-data objective by
\begin{equation}
 A=\bigl[\nabla^2R(\widehat b)+X^TDX\bigr]^{-1}.
 \label{eq:intro-smooth-inverse}
\end{equation}
Since $H_{-i}=A^{-1}-d_ix_ix_i^T$, the Sherman--Morrison identity
yields
\[
 x_i^TH_{-i}^{-1}x_i
 =\frac{x_i^TAx_i}{1-d_ix_i^TAx_i}.
\]
Thus the Newton approximation to the held-out prediction is
$x_i^T\widetilde b^{(-i)}
=t_i+s_ix_i^TAx_i/(1-d_ix_i^TAx_i)$.
This explains both the rational ALO correction and the expected
choice of $A$ for a $C^2$ penalty: it is the inverse of penalty
curvature plus fitted loss curvature
\cite{rad2018scalable,auddy2024approximate}.

\section{\texorpdfstring{$\sqrt n$}{sqrt(n)}-consistency of ALO and mean-field inference: The desired bounds}

For general penalties, we retain this prediction formula and seek a
matrix $A$ playing the same role. For an input matrix
$A\in\R^{p\times p}$ with nonzero denominators, define the ALO and
mean-field predictions by
\begin{align}
 t_i^{\rm ALO}(A)
 &=t_i+s_i\frac{x_i^TAx_i}{1-d_ix_i^TAx_i},&
 t_i^{\rm MF}(A)
 &=t_i+s_i\tr(A\Sigma).
 \label{eq:intro-predictions}
\end{align}
The mean-field prediction uses the common trace correction
$\tr(A\Sigma)$ in place of the
observation-specific rational weight. The corresponding criteria are
\begin{align}
 \alo(A)&=\frac1n\sum_{i=1}^n g(t_i^{\rm ALO}(A),y_i),&
 \operatorname{MF}(A)&=\frac1n\sum_{i=1}^n g(t_i^{\rm MF}(A),y_i).
 \label{eq:intro-criteria}
\end{align}
The choice of $A$ is central to both approximations.

In the proportional regime $p\le n/\delta$, the desired prediction
bounds are
\begin{equation}
 \frac1n\sum_{i=1}^n\E\left[
 |t_i^{\rm ALO}(A)-t_i^{\rm LOO}|^2
 +|t_i^{\rm MF}(A)-t_i^{\rm LOO}|^2\right]
 \le\frac{C(\mu,\delta)}n.
 \label{eq:intro-desired-prediction}
\end{equation}
Under the assumptions of \cref{sec:introduction} and the response
models in \cref{thm:intro-statistics}, Cauchy--Schwarz and prediction
moment bounds turn these into criterion bounds.
For the same choice of $A$ and every jointly measurable test function
satisfying $|g(t,y)-g(t',y)|\le|t-t'|(1+|t|+|t'|)$, we seek the
criterion bound
\begin{equation}
 \E|\alo(A)-\loo|+\E|\operatorname{MF}(A)-\loo|
 \le\frac{C'(\mu,\delta)}{\sqrt n}.
 \label{eq:intro-desired-criteria}
\end{equation}
These statements compare ALO and mean field with LOO; estimating
generalization error additionally requires a bound for LOO itself
\cite{rad2020error,xu2019consistent,zou2024theoretical}.

\section{Closed form for \texorpdfstring{$A$}{A} and results from Part~I: \texorpdfstring{$C^2$}{C2} penalty and group Lasso}

The matrix $A$ in Part~I~\cite{bellec2025simultaneous} describes how the penalized estimator responds
to changes in the design. Holding the responses fixed, Part~I~\cite{bellec2025simultaneous}
represents its coordinate derivatives, almost everywhere in $X$, as
\begin{equation}
 \frac{\partial\widehat b}{\partial x_{ij}}
 =-A\Big[L_{y_i}'(x_i^T\widehat b)e_j
                +d_i\widehat b_jx_i\Big].
 \label{eq:intro-sensitivity}
\end{equation}
Here $e_j$ is the $j$th coordinate vector.
Informally, $A$ converts the change in the first-order equations caused
by perturbing $x_{ij}$ into the resulting change in $\widehat b$.
Such derivative representations occur in \cite{bellec2021derivatives,bellec2022observable,bellec2020out_of_sample}.

For a smooth penalty the familiar expression is
\begin{equation}
 A=\bigl[\nabla^2R(\widehat b)+X^TDX\bigr]^{-1}.
 \label{eq:intro-C2-inverse}
\end{equation}
ALO has also been developed and analyzed for nonsmooth penalties
\cite{wang2018approximate,wang2018approximate_fast_tuning,rad2018scalable,lin2026accuracy}.
For group Lasso, let the groups $G_1,\ldots,G_K$ partition $\{1,\ldots,p\}$ and take
\begin{equation}
    \label{group-lasso-penalty}
 R(b)=\frac a2\norm b^2+\sum_{k=1}^K\lambda_k\norm{b_{G_k}},
 \qquad a>0,\quad\lambda_k>0.
\end{equation}
The covariance-scaled assumption \eqref{eq:intro-strong-convexity}
holds if $a\ge n\mu\|\Sigma\|_{op}$.
Let $\widehat S=\bigcup_{k:\widehat b_{G_k}\ne0}G_k$ be the active
support, and let $\Pi_{\widehat S}$ project onto its coordinate space.
Write $B_{\widehat S}$ for the Hessian of the ridge and active-group
penalties at $\widehat b$, extended by zero off $\widehat S$.
Its active blocks are
\[
 (B_{\widehat S})_{G_kG_k}
 =aI_{G_k}+\frac{\lambda_k}{\norm{\widehat b_{G_k}}}
 \left(I_{G_k}-\frac{\widehat b_{G_k}\widehat b_{G_k}^T}
                         {\norm{\widehat b_{G_k}}^2}\right).
\]
The corresponding matrix is
\begin{equation}
 A=\bigl[\Pi_{\widehat S}
          (B_{\widehat S}+X^TDX)\Pi_{\widehat S}\bigr]^\dagger.
 \label{eq:intro-group-inverse}
\end{equation}
This is the inverse on the active coefficient space, extended by zero;
it is zero when there are no active groups.

In both cases, the fit determines the penalty curvature and the
space on which it is inverted. Conditioning also fixes $D$, leaving
the explicit $X^TDX$ term as the remaining design dependence.
This structure is what the probabilistic comparison uses.

Earlier consistency analyses cover smooth penalties
\cite{rad2018scalable,xu2019consistent} and the elastic net
\cite{auddy2024approximate}. Part~I~\cite{bellec2025simultaneous}
established a probabilistic argument for
consistency in the group Lasso case \eqref{group-lasso-penalty}.
\citet{lin2026accuracy}
studied the error from restricting leave-one-out fits to the fitted active
space for blockwise penalty functions
of the form $R(b)=\sum_{k=1}^K\lambda_k r_k(b_{G_k})$ where $(G_1,...,G_K)$
are non-overlapping groups and $r_k$ are three times continuously
differentiable everywhere except at the origin
with subdifferential at the origin satisfying
$B(0,\rho_-)\subset \partial r_k(0) \subset B(0,\rho_+)$
for some $0<\rho_-<\rho_+<\infty$, under additional ridge,
block-size and regularity assumptions.

\section{The conditioning argument in Part~I}
\label{sec:part-i-conditioning}

Part~I~\cite{bellec2025simultaneous} first controls the mean-field prediction using the derivative
representation \eqref{eq:intro-sensitivity}, Gaussian integration by
parts and stability of the fit when one observation is removed:
\[
 \sum_{i=1}^n\E|t_i^{\rm MF}(A)-t_i^{\rm LOO}|^2\le C(\mu,\delta).
\]
This step needs derivative and operator-norm bounds for $A$.
In a single-index model, derivatives are taken on slices where the
external randomness and $X\beta$ are fixed, so that the responses
remain fixed. To transfer the conclusion to ALO, it remains to
compare the weights
\[
 W_i=\frac{x_i^TAx_i}{1-d_ix_i^TAx_i},
 \qquad \gamma=\tr(A\Sigma),
 \qquad
 t_i^{\rm ALO}(A)-t_i^{\rm MF}(A)=s_i(W_i-\gamma).
\]
The explicit inverse formulas
\eqref{eq:intro-C2-inverse} and \eqref{eq:intro-group-inverse}
make this comparison possible under the conditional law of the design.
In particular, fixing the group Lasso fit fixes its random active
support. The conditional law is generally not Gaussian; its
description is the next step.

For a single-index model $y_i=f_i(U,x_i^T\beta)$, first fix the
external randomness $U$, independent of $X$, and the index values
$a=X\beta$. Then condition on $\widehat b$, $t=X\widehat b$ and
the full penalty subgradient $q=-X^Ts$, including its inactive
coordinates. The scores $s$, curvature $D$, and responses are fixed.
For a smooth penalty this fixes $\nabla^2R(\widehat b)$; for group
Lasso it fixes $\widehat S$, $\Pi_{\widehat S}$ and $B_{\widehat S}$.
The KKT condition and uniqueness of the fit identify the fiber as
the affine space
\begin{equation}
 \mathcal F=\{X:X\beta=a,\ X\widehat b=t,\ X^Ts=-q\}.
 \label{eq:intro-conditioning-fiber}
\end{equation}
Fixing the entire subgradient makes the inactive KKT inequalities
fixed conditions on the conditioning values, so they impose no
additional restriction inside this fiber. For robust regression
with loss $\rho(y_i-x_i^Tb)$, one can instead fix the noise
$\varepsilon$, $h=\widehat b-\beta$, $u=Xh$ and
$q=X^T\rho'(\varepsilon-u)$; the fiber then has only the constraints
$Xh=u$ and $X^T\rho'(\varepsilon-u)=q$.

The Gaussian reference law on $\mathcal F$ has whitened
representation
$X\Sigma^{-1/2}=M+Q_\perp ZP_\perp$, where $Z$ has independent
standard normal entries and $M$ is its fixed mean.
Here $Q_\perp$ projects onto $s^\perp$ and $P_\perp$ onto
$\operatorname{span}\{\Sigma^{1/2}\beta,
\Sigma^{1/2}\widehat b\}^{\perp}$; zero or redundant directions are
omitted. Thus at most one left-hand and two right-hand directions
are removed. The actual conditional law $\pi$ includes the coarea
Jacobian of the nonlinear conditioning map. Under the smooth-loss
assumptions of Part~I~\cite{bellec2025simultaneous}, on regular fibers with nonempty active set,
it is the determinant tilt
\begin{equation}
 \frac{d\pi}{d\gamma_{\mathcal F}}(X)
 =
 \frac{\det(B_S+X_S^TDX_S)}
 {\E_{\gamma_{\mathcal F}}\det(B_S+X_S^TDX_S)}.
 \label{eq:intro-conditioning-law}
\end{equation}
Here $\gamma_{\mathcal F}$ denotes the Gaussian reference law,
and $X_S$ contains the columns of $X$ indexed by $S$.
In the smooth case, $S=[p]$ and $B_S=\nabla^2R(\widehat b)$.
For group Lasso, $S=\widehat S$; in this conditional-law formula
and the bounds below, $B_S$ denotes the $|S|\times|S|$ restriction
of the zero-extended penalty Hessian in
\eqref{eq:intro-group-inverse} to its active coordinates.
The conditional-law lemmas in
Part~I~\cite[Lemmas A.3 and B.1]{bellec2025simultaneous} also handle
the degenerate cases: empty active supports give zero weight error,
zero scores give an affine Gaussian conditional law, and collinear
directions are treated by a lower-dimensional calculation.
Thus the argument covers the degenerate events as well as the
regular fibers.

This describes the smooth-loss calculation. The formal argument
extends it to 1-Lipschitz scores using almost-everywhere derivatives
of the scores and the fitted KKT map, with the convention for $d_i$
given above. Thus bounded, 1-Lipschitz scores suffice for the
statistical results.

With $B_S\succeq n\mu\Sigma_{S,S}$, the gradient of this
log-determinant tilt in the free whitened coordinates has norm at
most $1/\sqrt\mu$. A weighted second-order Stein inequality therefore
controls, under $\pi$, the remainders
\[
 \Rem_i=x_i^TAx_i-\gamma(1-d_ix_i^TAx_i).
\]
Strong convexity also gives the pointwise denominator bound
\[
 0<(1-d_ix_i^TAx_i)^{-1}
 \le1+\frac{d_ix_i^T\Sigma^{-1}x_i}{n\mu}.
\]
Define the row-norm event
\[
 \Omega=\{\max_i\|\Sigma^{-1/2}x_i\|\le\sqrt p+\sqrt n\}.
\]
Let $I_\Omega$ denote its indicator.
Using $W_i-\gamma=\Rem_i/(1-d_ix_i^TAx_i)$, the conditional
comparison on $\Omega$ takes the form
\begin{equation}
 \E_\pi\left[I_\Omega\sum_{i=1}^n|W_i-\gamma|^2\right]
 \le C(\mu,\delta)\E_\pi\left[
 \left(1+\frac{\|X\Sigma^{-1/2}\|_{op}}{\sqrt n}\right)^4\right].
 \label{eq:intro-conditional-weights}
\end{equation}
For group Lasso the same estimate holds by restricting the inverse
to the now fixed active space.

Averaging over the conditioning values uses the tower property:
the right-hand side becomes a moment of the original Gaussian
design, bounded independently of $n,p$. Consequently
$\E[I_\Omega\sum_{i=1}^n|W_i-\gamma|^2]\le C(\mu,\delta)$.
Together with $|s_i|\le1$, the mean-field bound and
Cauchy--Schwarz, this gives the ALO prediction and criterion bounds
on $\Omega$. Part~II also controls the Gaussian tail contribution from
$\Omega^c$, giving the unconditional bounds
\eqref{eq:intro-desired-prediction}--\eqref{eq:intro-desired-criteria}.

\section{General convex penalty: Which matrix \texorpdfstring{$A$}{A}?}

The only obstacle to extending the argument of Part~I~\cite{bellec2025simultaneous} to general
finite-valued strongly convex penalties is therefore the absence
of an explicit formula for $A$ with this separation of
dependences. We seek the same structure for a
general penalty: a response space and penalty curvature determined
by the fitted quantities, with $X$ appearing explicitly inside
the inverse or pseudoinverse on that space.

The two contributions of Part~II are:
\begin{enumerate}
\item Explicit expressions for $A$ under general penalties.
We construct explicit matrices for general finite-valued strongly
convex penalties, recovering the familiar $C^2$ and group Lasso
matrices in those cases. Several constructions are possible; this
introduction presents only $A_{\rm der}$, based on ordinary derivatives
of the penalty response.
\item A formal proof of the statistical guarantees.
For this new expression of $A$, we prove in Lean~4 that the
probabilistic argument of Part~I~\cite{bellec2025simultaneous} goes through and establishes the
success of ALO under the same model and score assumptions. We also
formally prove the mean-field guarantees for the same choice of $A$.
The distinction is that the mean-field argument does not require an
explicit expression for $A$, whereas the ALO argument uses it.
\end{enumerate}
The later sections give the proofs and describe the other constructions.

\section{From derivatives of the penalty Fenchel conjugate
to the proposed \texorpdfstring{$A_{\rm der}$}{A(der)}}
Define the \emph{penalty response} by
\[
    G_R(v)=\argmin_{b \in \R^p}\{R(b)-v^Tb\}.
\]
The minimum value is $-R^*(v)$, where $R^*$ is the Fenchel
conjugate of $R$, and $G_R(v)=\nabla R^*(v)$.
Strong convexity of $R$ makes $G_R$ Lipschitz. This response
asks how the coefficient preferred by the penalty changes under a
linear force $v$. In particular $G_R(q)=\widehat b$. Choose a linear
space $F$ of largest dimension on which this response has an ordinary
first-order expansion at $q$:
\begin{equation}
 G_R(q+z)=\widehat b+J_Fz+o(\norm z),\qquad z\in F,\quad z\to0.
 \label{eq:intro-response}
\end{equation}
If the usual Jacobian exists, take $F=\R^p$. Otherwise the restriction
allows the construction to retain directions with a linear response.
A fixed measurable rule resolves ties; the zero space is always
available. Let $P$ be orthogonal projection onto $F$ and put $J=J_FP$.
Our choice is
\begin{equation}
 \boxed{\displaystyle
 A_{\rm der}=J\bigl(PJ+J^TX^TDXJ\bigr)^\dagger J^T,}
 \label{eq:intro-ader}
\end{equation}
where $\dagger$ denotes the Moore--Penrose pseudoinverse. To display
the design dependence in the same form as the preceding examples, put
$T=\ran J$ and let $\Pi$ project onto $T$. The response determines a
symmetric positive definite penalty curvature $H_{\rm der}:T\to T$,
characterized by $H_{\rm der}Jz=\Pi Pz$. Writing
$\overline H_{\rm der}$ for its extension by zero off $T$, we have
\begin{equation}
 \boxed{\displaystyle
 A_{\rm der}
 =\bigl[\Pi(\overline H_{\rm der}+X^TDX)\Pi\bigr]^\dagger.}
 \label{eq:intro-ader-inverse}
\end{equation}
Thus we invert penalty curvature plus loss curvature on the response
space, and extend by zero on its orthogonal complement. Both $\Pi$ and
$\overline H_{\rm der}$ depend only on $R,q$. Once the conditioning
fixes $q$ and $D$, the remaining randomness is the $X$ appearing
explicitly inside this pseudoinverse, just as for the smooth and
group Lasso matrices. The equivalent formulas are proved in
\cref{prop:ordinary-penalty-selection}; see
\eqref{eq:ordinary-penalty-inverse}.

The matrix is Borel measurable and positive semidefinite. Constructing
$A_{\rm der}$ requires neither the design covariance nor an unknown
signal parameter. It recovers both matrices displayed above exactly.

\begin{proposition}[Recovery of the smooth and group Lasso expressions]
\label{prop:intro-recovery}
Let $R$ be finite valued and strongly convex, and construct
$A_{\rm der}$ by \eqref{eq:intro-ader} at the fit \eqref{eq:intro-fit}.
\begin{enumerate}
\item If $R$ is twice continuously differentiable, then
\[
 A_{\rm der}=\bigl[\nabla^2R(\widehat b)+X^TDX\bigr]^{-1}.
\]
\item For the ridge-augmented group Lasso penalty above, with its active
support and curvature notation,
\[
 A_{\rm der}
 =\bigl[\Pi_{\widehat S}
       (B_{\widehat S}+X^TDX)\Pi_{\widehat S}\bigr]^\dagger.
\]
This includes empty active sets and inactive groups at the boundary of
their KKT conditions; no strict complementarity is required.
\end{enumerate}
Both identities hold pointwise, rather than only almost surely or
asymptotically.
\end{proposition}

The proofs and their Lean~4 formalizations are given in
\cref{sec:ordinary-penalty}, including
\eqref{eq:ordinary-penalty-group-recovery}.

\section{Results for general convex penalties and the data-driven \texorpdfstring{$A_{\rm der}$}{A(der)}}

We use $A=A_{\rm der}$ in the theorem below. For this choice, ALO
uses only the observed data, loss and penalty; the mean-field correction
additionally uses $\Sigma$.

\begin{theorem}[ALO and mean field for a general penalty]
\label{thm:intro-statistics}
Let $n\ge1$, $\mu,\delta>0$, and $p\le n/\delta$. Suppose the rows of
$X$ are independent $N(0,\Sigma)$ vectors with $\Sigma\succ0$.
Let $R$ be a deterministic penalty with fixed tuning parameters,
possibly depending on $n,p$. Assume $R$ is finite valued, minimized
at zero, and such that
$R(b)-n\mu b^T\Sigma b/2$ is convex. Suppose the losses are convex and
differentiable, with jointly measurable scores $\psi_y=L_y'$ satisfying
$|\psi_y|\le1$ and $\operatorname{Lip}(\psi_y)\le1$.
Allow any single-index response model
\[
 y_i=f_i(U,x_i^T\beta),
\]
where $\beta$ is fixed, the links $f_i$ are measurable, and $U$ is
independent of $X$ with an arbitrary probability law. This includes
linear regression $y_i=x_i^T\beta+\eps_i$ with an arbitrary joint noise
law independent of $X$.

For $A=A_{\rm der}$ the denominators in \eqref{eq:intro-predictions}
are strictly positive, and
\[
 0\preceq A_{\rm der}\preceq\frac{\Sigma^{-1}}{n\mu}.
\]
There are finite constants $C(\mu,\delta)$ and $C'(\mu,\delta)$ such
that, for every jointly measurable test function satisfying
$|g(t,y)-g(t',y)|\le|t-t'|(1+|t|+|t'|)$,
\begin{align}
    \frac1n\sum_{i=1}^n\E\Big[
     |t_i^{\rm ALO}(A_{\rm der})-t_i^{\rm LOO}|^2
     +|t_i^{\rm MF}(A_{\rm der})-t_i^{\rm LOO}|^2\Big]
 &\le\frac{C(\mu,\delta)}{n},\label{eq:intro-prediction-rate}\\
 \E\big|\alo(A_{\rm der})-\loo\big|
 +\E\big|\operatorname{MF}(A_{\rm der})-\loo\big|
 &\le\frac{C'(\mu,\delta)}{\sqrt n}.\label{eq:intro-criterion-rate}
\end{align}
These are integrable errors under the full original noise/design law;
no event restricting the observed data is imposed.
\end{theorem}

Examples of admissible test functions are the absolute error
$g(t,y)=|y-t|$ in regression and the logistic loss
$g(t,y)=\log(1+\exp(-yt))$ in classification with $y\in\{-1,1\}$.
Both are 1-Lipschitz in $t$, uniformly in $y$, and therefore satisfy
the test-function condition.

Thus the general-penalty construction preserves both ALO accuracy and
the mean-field approximation, rather than providing only a sensitivity
formula. For ALO, the proof uses the explicit inverse structure, a finite-rank
comparison with the derivative operator, and the probabilistic estimates
of Part~I~\cite{bellec2025simultaneous}. The mean-field proof uses derivative and trace controls,
without requiring an explicit formula for $A$.

\section{Lean~4 proofs}

We identify each recovery identity in \cref{prop:intro-recovery}
and each statistical conclusion in \cref{thm:intro-statistics}
by its Lean declaration below, together with the supporting
construction and conditioning statements from Sections~6 and~4.
All names belong to the namespace \lean{ALO}; the source files are
in \lean{anc/lean_alo/ALO/Estimator}.

\paragraph{The recovery proposition.}
The smooth-penalty identity in \cref{prop:intro-recovery}(1) is
\lean{ordinaryPenaltyALOOperator_eq_smooth_inverse}, in
\lean{OrdinaryPenaltyAlgebra.lean}.
The group Lasso identity in \cref{prop:intro-recovery}(2) is
\lean{ordinaryPenaltyALOOperator_eq_groupLassoEmpiricalALOOperator},
in \lean{OrdinaryPenaltyGroupRecovery.lean}.
The underlying restricted inverse is identified by
\lean{ordinaryPenaltyALOOperator_eq_groupLasso_reducedInverse}.
These are pointwise identities for positive group weights; their
statements allow an empty active space and inactive groups at
the boundary of the subgradient conditions.

\paragraph{The statistical theorem.}
For \cref{thm:intro-statistics}, the single-index ALO prediction and
criterion bounds are the two conclusions of
\lean{singleIndex_ordinaryPenalty_ALO_unconditional};
the corresponding mean-field bounds are
\lean{singleIndex_ordinaryPenalty_meanField_unconditional}.
Both are in \lean{OrdinaryPenaltySingleIndex.lean}.
Their regression specializations are
\lean{regression_ordinaryPenalty_ALO_unconditional} and
\lean{regression_ordinaryPenalty_meanField_unconditional},
in \lean{OrdinaryPenaltyRegression.lean}.
Each statement includes integrability of the summed squared
prediction error and the absolute criterion error under the full
independent-noise/Gaussian-design law. Dividing the summed prediction
bounds by $n$ and adding the ALO and mean-field bounds gives
\eqref{eq:intro-prediction-rate}; adding their criterion bounds gives
\eqref{eq:intro-criterion-rate}. Thus these declarations also establish
the desired rates of Section~2 for $A=A_{\rm der}$.

Equation~\eqref{eq:intro-criterion-rate} also has the separate declaration
\lean{manuscript_theorem_7_1_equation_7_2}, in
\lean{ManuscriptCriterion.lean}. Its statement specifies the original
full-data and deleted minimization problems, defines $A_{\rm der}$ by
the projection--pseudoinverse formula, and writes each criterion as
its own average. It derives the combined bound from the two declarations
above, with a constant depending only on $\mu,\delta$.

The pointwise assertions of \cref{thm:intro-statistics} are supplied
by \lean{ordinaryPenaltyALOOperator_control},
\lean{ordinaryPenaltyALOOperator_spec},
\lean{norm_whiten_ordinaryPenaltyALOOperator_le} and
\lean{ordinaryPenaltyALOWeight_bound}: these give symmetry,
positive semidefiniteness, the covariance-relative norm bound,
strict positivity of the rational denominators and the weight bound.
The statistical proofs use whitened coordinates.
\lean{unwhiten_covarianceOrdinaryPenaltySignalOperator}, in
\lean{OrdinaryPenaltyCovarianceObserved.lean}, identifies their
operator with the original observed-data selection;
\lean{singleIndex_observedOrdinaryPenaltyOperator_eq_original_loss}
identifies its inputs with $q=-X^Ts$ and $X^TDX$ for the original
losses. These identities ensure that the endpoints concern the
matrix defined in \eqref{eq:intro-ader}.

\paragraph{The construction in Section~6.}
In \lean{OrdinaryPenaltyResponse.lean},
\lean{exists_ordinaryPenaltySelection} proves existence of a Borel
selection, and \lean{ordinaryPenaltySelection_spec} proves the
ordinary first-order expansion \eqref{eq:intro-response} and maximal
dimension of $F$.
\lean{measurable_ordinaryPenaltyALOOperator} proves measurability of
the resulting matrix.
In \lean{OrdinaryPenaltyAlgebra.lean},
\lean{ordinaryPenaltySelection_curvature} proves existence and
uniqueness of the symmetric coercive curvature $H_{\rm der}$;
\lean{ordinaryPenaltyALOOperator_eq_responseGram} proves
\eqref{eq:intro-ader}, and
\lean{ordinaryPenaltyALOOperator_eq_symmetricPseudoInverse} proves
\eqref{eq:intro-ader-inverse}.
The Lipschitz property of $G_R$ is
\lean{lipschitzWith_penaltyResponse}, in \lean{PenaltyResponse.lean}.

\paragraph{The conditioning estimates reviewed in Section~4.}
The determinant description of the actual conditional law on regular
fibers is \lean{ae_isSignalKKTDeterminantLaw_condDistrib}, in
\lean{SignalSliceGlobalDisintegration.lean}.
\lean{map_signalSliceKKTData_signalKKTDeterminantConditionalKernel}
states that this kernel fixes the fitted tuple, and
\lean{lintegral_joint_signalKKTDeterminantConditionalKernel} gives
the disintegration identity used to average conditional estimates.
The conditional remainder bound is
\lean{ae_signalSliceRemainder_bound_condDistrib}, in
\lean{SignalSliceConditionalRemainder.lean}.
After the denominator estimate and averaging,
\lean{constructedSignal_weight_uniform_bound} gives the summed
ALO-weight versus trace comparison on the row-norm event.
The mean-field-to-LOO prediction estimate is
\lean{signalAdaptiveConstructed_meanField_error_bound}, in
\lean{SignalAdaptiveMeanFieldModel.lean}.
These component statements are formulated for the derivative
operator; the finite-rank comparison
\lean{ae_covarianceOrdinaryPenaltySignalOperator_rank} transfers
their consequences to $A_{\rm der}$ in the statistical endpoints
above. The supplement gives the full finite-sample bounds,
constants and comparison statements in
\cref{sec:ordinary-penalty,sec:unconditional}.

\section*{Lean code availability}

The Lean~4 proofs are available in the ancillary files of this
paper's arXiv submission, in \lean{anc/lean_alo/}. The project
includes the source files and pinned dependency versions needed
to build and check the formalization.

\bibliographystyle{plainnat}
\bibliography{../../bibliography/db}

\section*{Supplement}
\renewcommand{\appendixname}{Supplement}
\appendix

\section{The canonical operator and its properties}
\label{sec:operator}

All vector spaces are real and finite dimensional. Rows of
$X\in\R^{n\times p}$ are denoted by $x_i^T$.
Fix an invertible matrix $S$ with $\Sigma=S^TS\succ0$, and write
\begin{equation}
 Z=XS^{-1},\qquad z_i=S^{-T}x_i,\qquad
 \theta=Sb,\qquad \widetilde R(\theta)=R(S^{-1}\theta).
 \label{eq:whitening}
\end{equation}
Assume that $R$ is finite valued and strongly convex, and that
$\widetilde R$ is $\alpha$-strongly convex for some $\alpha>0$.
For a fixed $c\in\R^p$ and convex differentiable functions
$\ell_i:\R\to\R$ whose derivatives are 1-Lipschitz, define
\begin{equation}
 \widehat\theta(Z)=\argmin_{\theta\in\R^p}
 \left\{\sum_{i=1}^n\ell_i\bigl(z_i^T(\theta-Sc)\bigr)
       +\widetilde R(\theta)\right\},\qquad
 \widehat b=S^{-1}\widehat\theta.
 \label{eq:estimator}
\end{equation}
The coefficient is always optimized over the entire space $\R^p$.

The parameter with respect to which we differentiate is the design on
an affine space $\mathcal A=Z_0+\mathcal E$. The two choices used below are
\begin{equation}
 \mathcal E=\R^{n\times p},
 \qquad\text{or}\qquad
 \mathcal E=\{E\in\R^{n\times p}:Es=0\}
 \quad\text{for a fixed }s\in\R^p.
 \label{eq:design-directions}
\end{equation}
In the second case, $Zs=Z_0s$ is fixed. If $s=0$, this is the full
design space. In defining the derivative on $\mathcal A$, the loss
functions $\ell_i$ and the offset $c$ are held fixed.

Set
\begin{equation}
 h=\widehat\theta-Sc,\quad u=Zh,\quad
 v_i=-\ell_i'(u_i),\quad d_i=(\ell_i')'(u_i),\quad
 D=\diag(d_1,\ldots,d_n).
 \label{eq:kkt-variables}
\end{equation}
We set $(\ell_i')'(t)=0$ when this derivative does not exist, as in
Lean's \lean{deriv}. Convexity and the Lipschitz assumption give
$0\le d_i\le1$. The derivative of a Lipschitz score composed with a
fitted predictor obeys the required chain rule almost everywhere,
including this convention at score kinks.

\begin{definition}[Canonical factor on the chosen design space]
\label{def:operator}
Let $B:\mathcal E\to\R^p$ be the derivative of
$Z\mapsto\widehat\theta(Z)$ along $\mathcal A$, with value zero where
the derivative does not exist. Define the forcing map
\begin{equation}
 \mathcal L(E)=E^Tv-Z^TDEh,\qquad E\in\mathcal E.
 \label{eq:forcing}
\end{equation}
The matrix $\widetilde A$ is the unique solution, when one exists, of
\begin{equation}
 \widetilde A^T=\widetilde A,\qquad
 \widetilde A\mathcal L=B,\qquad
 \ran\widetilde A=\ran B.
 \label{eq:factor-definition}
\end{equation}
Set $\widetilde A=0$ if there is no solution. The operator in the
original coefficient coordinates is
\begin{equation}
 A=S^{-1}\widetilde A S^{-T}.
 \label{eq:original-operator}
\end{equation}
\end{definition}

The range condition in \eqref{eq:factor-definition} is part of the
definition: the derivative equation alone need not identify a matrix.
This operator depends on the allowed design perturbations, as well as
on the loss and penalty. In particular, the operator on a fixed-signal
slice is constructed from the derivative on that slice.

\begin{remark}[Oracle dependence and the data-driven replacement]
\label{rem:oracle}
For the regression application below the offset is the unknown true
coefficient, and for the single-index application the allowed design
directions depend on the unknown signal direction. Consequently
\cref{def:operator} does not establish an operator determined by
$(y,X)$ alone. Borel measurability for each fixed model parameter does
not remove this dependence.
The main statistical theorems in \cref{sec:unconditional} use the
data-driven replacement $A_{\rm data}$ of \cref{sec:data-operator}.
It is defined in observed coordinates and requires no knowledge of
$\Sigma$. The canonical operator above supplies the derivative and
statistical comparisons used in its proof. Independently, the prescribed
group Lasso matrix in \cref{def:group-inverse} has the unconditional
bounds of \cref{cor:prescribed-group}. Exact equality of $A_{\rm data}$
with the prescribed smooth and positive-weight group matrices is proved in
\cref{thm:data-smooth-recovery,thm:data-group-recovery}.
\end{remark}

\begin{theorem}[Existence, uniqueness, and restricted inverse]
\label{thm:operator-existence}
Under the preceding assumptions, the solution of
\eqref{eq:factor-definition}, whenever it exists, is unique.
The map $Z\mapsto\widetilde A(Z)$ is Borel measurable on each fixed
affine design space, and the solution exists for Lebesgue-almost every
$Z$ in that space.

More precisely, at almost every such $Z$, put $T=\ran B$ and let
$P_T$ be orthogonal projection onto $T$. The KKT subgradient
$q(Z)=Z^Tv(Z)\in\partial\widetilde R(\widehat\theta(Z))$
is differentiable along $\mathcal A$. There is a unique symmetric
operator $H:T\to T$ satisfying
\begin{equation}
 H\,B(E)=P_T\,D_{\mathcal E}q(Z)[E]\quad(E\in\mathcal E),
 \qquad H\succeq\alpha I_T.
 \label{eq:tangent-curvature}
\end{equation}
If $U_T$ has orthonormal columns $e_1,\ldots,e_r$ spanning $T$, then
\begin{equation}
 \widetilde A
 =U_T\bigl(H_T+U_T^TZ^TDZU_T\bigr)^{-1}U_T^T,
 \qquad H_T=(e_j^THe_k)_{j,k=1}^r.
 \label{eq:restricted-inverse}
\end{equation}
Here $H_T$ denotes the matrix of $H$ in this basis; the expression is
zero when $T=\{0\}$. Thus no Hessian of the penalty is needed.
\end{theorem}

\begin{proof}
Strong convexity applied to the differentiable coefficient/subgradient
pair yields the symmetric operator in \eqref{eq:tangent-curvature}.
The differentiated KKT equation is
$D_{\mathcal E}q=\mathcal L-Z^TDZB$.
Projection onto $T$ gives \eqref{eq:restricted-inverse} and
$B=\widetilde A\mathcal L$. Its range is $T$.
Two symmetric factors satisfying the range condition in
\eqref{eq:factor-definition} agree. A Borel selection of this unique
factor gives the specified extension at exceptional designs.
\end{proof}

\paragraph{Lean formalization.}
The definition is \lean{measurableEmpiricalALOOperator}, using
\lean{symmetricFactor}. The relevant results include
\lean{IsSymmetricFactor.unique},
\lean{measurable_symmetricFactor},
\lean{existsUnique_curvature_of_strongConvex_local}, and
\lean{derivative_eq_reducedInverse}.
Almost-everywhere existence with Lipschitz scores is
\lean{ae_measurableEmpiricalALOOperator_spec_of_lipschitz_score}.
The slice construction and its specification are
\lean{signalSliceOperator} and \lean{ae_signalSliceOperator_spec}.
These declarations are in the \lean{ALO} namespace.

\begin{theorem}[Derivative, energy, and denominator bounds]
\label{thm:operator-properties}
Under the assumptions of \cref{thm:operator-existence}, the following
hold almost everywhere on the chosen design space:
\begin{align}
 D_{\mathcal E}\widehat\theta(Z)[E]
   &=\widetilde A(E^Tv-Z^TDEh),\qquad E\in\mathcal E,
   \label{eq:derivative}\\
 0\preceq\widetilde A&\preceq\alpha^{-1}I_p,\qquad
 \widetilde A Z^TDZ\widetilde A\preceq\widetilde A,\qquad
 \ran\widetilde A=\ran D_{\mathcal E}\widehat\theta(Z).
 \label{eq:energy}
\end{align}
Consequently, in original coordinates,
\begin{equation}
 0\preceq A\preceq\alpha^{-1}\Sigma^{-1},\qquad
 AX^TDXA\preceq A,\qquad
 0\le\tr(A\Sigma)=\tr(\widetilde A)\le p/\alpha.
 \label{eq:covariance-energy}
\end{equation}
For every row, let $a_i=x_i^TAx_i=z_i^T\widetilde A z_i$. Then
\begin{align}
 1-d_i a_i&\ge
       \frac{\alpha}{\alpha+d_i\norm{z_i}^2}>0,
       \label{eq:denominator}\\
 (1-d_i a_i)^{-1}&\le1+\frac{d_i\norm{z_i}^2}{\alpha},
 &0\le W_i(A):=\frac{a_i}{1-d_i a_i}
       &\le\frac{\norm{z_i}^2}{\alpha}.
       \label{eq:rational-bound}
\end{align}
These bounds hold without restricting to a row-norm event.
\end{theorem}

\paragraph{Lean formalization.}
The derivative, range, positivity, norm, energy, and denominator
conclusions are the fields of \lean{IsEmpiricalALO}, established by
the Lipschitz-score specification cited above.
The rational bound is \lean{operator_rational_weight_bound}.
For the random single-index slice, measurability and integration over
the signal coordinates are established separately; see
\lean{measurable_signalSliceOperator},
\lean{ae_signalAdaptive_regular}, and
\lean{ae_signalAdaptive_denominator}.
Thus the slice assertions are obtained under the original Gaussian law,
rather than by evaluating a fixed-loss null-set assertion at random losses.

\section{A data-driven operator in pseudoinverse form}
\label{sec:data-operator}

This construction uses the original coefficient coordinates. Let
$R:\R^p\to\R$ be finite valued and $a$-strongly convex for some $a>0$.
For convex differentiable observed losses $L_{y_i}$ with 1-Lipschitz
derivatives, set
\begin{align}
 \widehat b&=\argmin_b
   \left\{\sum_{i=1}^n L_{y_i}(x_i^Tb)+R(b)\right\},\\
 s_i&=L_{y_i}'(x_i^T\widehat b),\qquad
 q=-X^Ts\in\partial R(\widehat b),\\
 M&=X^TDX,\qquad D=\diag(d_1,\ldots,d_n).
 \label{eq:data-inputs}
\end{align}
Here $d_i=(L_{y_i}')'(x_i^T\widehat b)$, with value zero where this
derivative does not exist. Thus $0\le d_i\le1$. The inputs $R,q,M$
are determined by the penalty, loss and observed data.

\begin{definition}[Selected linear penalty response]
\label{def:data-response}
Define the penalty response
\begin{equation}
 G_R(u)=\argmin_{b\in\R^p}\{R(b)-u^Tb\}.
 \label{eq:penalty-response}
\end{equation}
For each $q$, select a locally uniform subsequential limit $g_q$ of
\begin{equation}
 g_{q,m}(z)=m\bigl[G_R(q+z/m)-G_R(q)\bigr],
 \qquad m=1,2,\ldots.
 \label{eq:data-blowup}
\end{equation}
Select a linear subspace $F_q\subseteq\R^p$ of largest dimension on
which $g_q$ is linear, and put
\begin{equation}
 P=P_{F_q},\qquad Jz=g_q(Pz),\qquad T_q=\ran J.
 \label{eq:data-response-matrices}
\end{equation}
Both selections are made by fixed Borel rules. Fix these rules for
the penalty $R$ before evaluating them at the observed $q$.
\end{definition}

The two subspaces have different roles: $F_q$ consists of subgradient
perturbations, whereas $T_q$ consists of coefficient responses. Write
$\Pi=P_{T_q}$ for the latter orthogonal projector; in general $P\ne\Pi$.
The modulus $a$ is used to prove the properties below, and is not an
input to the selection rule.

\begin{theorem}[Measurable response and induced curvature]
\label{thm:data-response}
The selections in \cref{def:data-response} exist, and $q\mapsto(P,J)$
is Borel measurable. They satisfy
\begin{equation}
 J=JP,\qquad PJ=J^TP,\qquad
 z^TPJz\ge a\norm{Jz}^2\quad(z\in\R^p).
 \label{eq:data-response-energy}
\end{equation}
There is a unique symmetric operator $H_q:T_q\to T_q$ such that
\begin{equation}
 H_qJz=\Pi Pz\quad(z\in\R^p),\qquad H_q\succeq aI_{T_q}.
 \label{eq:data-curvature}
\end{equation}
\end{theorem}

\begin{proof}
The response $G_R$ is globally $a^{-1}$-Lipschitz and
$G_R(q)=\widehat b$. The anchored maps in \eqref{eq:data-blowup}
form a relatively compact family for locally uniform convergence.
A Borel cluster-point selection gives $g_q$. Linearity on a subspace
is encoded by its orthogonal projector and the matrix $J$.
These matrix pairs form compact fibers. Minimizing minus the projector
trace maximizes the dimension; fixed separating coordinates resolve
ties measurably. Strong monotonicity and cyclic monotonicity of the
penalty subgradient pass to the limit and yield
\eqref{eq:data-response-energy}. The resulting linear response
relation gives the symmetric coercive curvature
\eqref{eq:data-curvature} on $\ran J$.
\end{proof}

\begin{definition}[Data-driven ALO operator]
\label{def:data-operator}
Extend $H_q$ by zero on $T_q^\perp$, writing the extension as
$\overline H_q$. Let ${}^\dagger$ denote the Moore--Penrose
pseudoinverse. Define
\begin{equation}
 \boxed{\displaystyle
 A_{\rm data}
   =\bigl[\Pi(\overline H_q+X^TDX)\Pi\bigr]^\dagger.}
 \label{eq:data-projection-pseudoinverse}
\end{equation}
Equivalently, if $U$ has orthonormal columns spanning $T_q$,
\begin{equation}
 A_{\rm data}
   =U\bigl[U^T(\overline H_q+X^TDX)U\bigr]^{-1}U^T.
 \label{eq:data-restricted-inverse}
\end{equation}
Set $A_{\rm data}=0$ when $T_q=\{0\}$.
\end{definition}

The compressed matrix is positive definite on $T_q$ and zero on
$T_q^\perp$, which explains the pseudoinverse in
\eqref{eq:data-projection-pseudoinverse}. In precisely the same notation,
the prescribed group Lasso matrix of \cref{def:group-inverse} is
\begin{equation}
 \boxed{\displaystyle
 Q_{\rm grp}
   =\bigl[P_{T_{\rm grp}}(H_{\rm grp}+X^TDX)P_{T_{\rm grp}}\bigr]^\dagger.}
 \label{eq:group-projection-pseudoinverse}
\end{equation}
Thus both constructions invert the loss curvature plus a penalty
curvature on a selected coefficient subspace and extend the result by
zero. This common form does not by itself identify $T_q,H_q$
with $T_{\rm grp},H_{\rm grp}$.

\begin{proposition}[Formula using only the selected response matrices]
\label{prop:data-factor-pseudoinverse}
The same operator can be written without explicitly constructing $H_q$:
\begin{equation}
 \boxed{\displaystyle
 A_{\rm data}
   =J\bigl(PJ+J^TX^TDXJ\bigr)^\dagger J^T.}
 \label{eq:data-factor-pseudoinverse}
\end{equation}
\end{proposition}

\begin{proof}
Let $N=\overline H_q+M$ and $C=J^TNJ$.
By \eqref{eq:data-curvature} and $PJ=J^TP$,
$J^T\overline H_qJ=J^TP=PJ$, so
$C=PJ+J^TMJ$. Since $N$ is positive definite on $\ran J$,
$C\succeq0$ and $\ker C=\ker J$. Consequently
$JC^\dagger C=J$. The symmetric matrix $B=JC^\dagger J^T$
has range contained in $T_q$, and
\[
 B(\Pi N\Pi)J=JC^\dagger J^TNJ=J.
\]
As $J$ maps onto $T_q$, $B$ is the inverse of the compressed curvature
on $T_q$ and is zero on $T_q^\perp$. This is
\eqref{eq:data-projection-pseudoinverse}.
\end{proof}

\paragraph{Lean formalization.}
\lean{dataDrivenResponse} and \lean{dataDrivenALOOperator}
implement the selection and operator using only $R,q,M$.
In Lean the operator is defined as the unique symmetric factor of $J$
through $P+MJ$, with range $\ran J$:
\begin{equation}
 A_{\rm data}(P+MJ)=J,\qquad
 A_{\rm data}^T=A_{\rm data},\qquad \ran A_{\rm data}=\ran J.
 \label{eq:data-factor-definition}
\end{equation}
The selection and response properties are established by
\lean{exists_penaltyResponseSelection},
\lean{dataDrivenResponse_spec}, and
\lean{dataDrivenResponse_linear_spec}.
The equality with \eqref{eq:data-projection-pseudoinverse} follows from
\lean{linearResponseOperator_eq_symmetricPseudoInverse}; the four
Moore--Penrose equations are checked in
\lean{restrictedInverse_moorePenrose}.
The explicit factor identity \eqref{eq:data-factor-pseudoinverse} is
proved above by finite-dimensional algebra; it has not yet been packaged
as a separate Lean equality. Exact recovery of the prescribed smooth
and positive-weight group Lasso matrices is checked in
\cref{thm:data-smooth-recovery,thm:data-group-recovery}.

\begin{theorem}[Covariance bounds and admissible ALO weights]
\label{thm:data-weights}
The map $(q,M)\mapsto A_{\rm data}$ is Borel measurable.
Suppose additionally that $R\circ S^{-1}$ is $\alpha$-strongly convex,
where $\Sigma=S^TS\succ0$. Then
\begin{align}
 0\preceq A_{\rm data}&\preceq\alpha^{-1}\Sigma^{-1},\\
 1-d_ix_i^TA_{\rm data}x_i&>0,\qquad
 (1-d_ix_i^TA_{\rm data}x_i)^{-1}
 \le1+\frac{d_i\norm{S^{-T}x_i}^2}{\alpha},\\
 0\le W_i(A_{\rm data})
 &=\frac{x_i^TA_{\rm data}x_i}{1-d_ix_i^TA_{\rm data}x_i}
 \le\frac{\norm{S^{-T}x_i}^2}{\alpha}.
 \label{eq:data-weight-bound}
\end{align}
These assertions hold at every observed dual point, with no row event.
\end{theorem}

\paragraph{Lean formalization and statistical scope.}
These properties are
\lean{measurable_dataDrivenALOOperator},
\lean{dataDrivenALOOperator_spec},
\lean{norm_whiten_dataDrivenALOOperator_le},
\lean{dataDrivenALOOperator_denominator}, and
\lean{dataDrivenALOWeight_bound}.
The fully instantiated standard-Gaussian signal-adaptive endpoints are
\lean{dataDrivenSignal_ALO_test_unconditional} and
\lean{dataDrivenSignal_meanField_unconditional}, both without $I_\Omega$.
The fully instantiated general-covariance endpoints are
\lean{dataDrivenCovariance_ALO_unconditional} and
\lean{dataDrivenCovariance_meanField_unconditional}. They use the
original-coordinate selection under the full independent-noise/design law,
with no derivative, rank, comparison, or moment premises left to supply.
Named specializations include
\lean{regression_dataDriven_ALO_unconditional},
\lean{singleIndex_dataDriven_ALO_unconditional}, and
\lean{singleIndex_groupLasso_dataDriven_ALO_unconditional}, with
corresponding mean-field endpoints.
The identity
\lean{singleIndex_observedDataDrivenOperator_eq_original_loss}
expresses the operator directly through the realized response scores and
the ordinary empirical fit. These data-driven endpoints are stated
in \cref{thm:linear,thm:single,cor:data-group}, including their
mean-field conclusions.

Computing $A_{\rm data}$ and its rational weights uses neither $\Sigma$
nor a true coefficient or signal direction. The separate mean-field
correction $\tr(A_{\rm data}\Sigma)$ still involves $\Sigma$.
The general construction is an existence-level definition; the
pseudoinverse formula does not supply a finite algorithm for finding
the selected blowup and maximal linearity subspace of an arbitrary penalty.

\section{Ordinary derivatives of the penalty response}
\label{sec:ordinary-penalty}

There is a second construction using only the original-coordinate
inputs $R,q,M$ of \cref{sec:data-operator}, with no selected
subsequential limit. It differentiates the penalty response itself
on a maximal subspace. This changes the selection rule for a general
strongly convex penalty, but preserves the statistical guarantees and
exact recovery of the prescribed smooth and group Lasso matrices.

\begin{definition}[Maximal ordinary penalty-response derivative]
\label{def:ordinary-penalty-response}
Let $G_R$ be as in \eqref{eq:penalty-response}. At each $q$, call a
linear subspace $F\subseteq\mathbb R^p$ admissible if there is a
linear map $J_F:F\to\mathbb R^p$ such that
\begin{equation}
 G_R(q+z)=G_R(q)+J_Fz+o(\norm z),\qquad z\in F,\quad z\to0.
 \label{eq:ordinary-penalty-expansion}
\end{equation}
This is ordinary Fr\'echet differentiability on $F$, along the full
limit $z\to0$. Choose an admissible $F$ of largest dimension, resolving
ties by a fixed Borel rule depending only on $R$, and put
\[
 P=P_F,\qquad J=J_FP,\qquad T=\ran J,\qquad \Pi=P_T.
\]
For $M=X^TDX$, define
\begin{equation}
 \boxed{\displaystyle
 A_{\rm der}=J\bigl(PJ+J^TMJ\bigr)^\dagger J^T.}
 \label{eq:ordinary-penalty-factor}
\end{equation}
The zero space is admissible, and its derivative and operator are zero.
\end{definition}

\begin{proposition}[Selection and compressed curvature]
\label{prop:ordinary-penalty-selection}
The selection can be made Borel in $q$. It satisfies
\[
 J=JP,\qquad PJ=J^TP,\qquad
 z^TPJz\ge a\norm{Jz}^2.
\]
There is a unique symmetric $H_{\rm der}:T\to T$ such that
\begin{equation}
 H_{\rm der}Jz=\Pi Pz,\qquad H_{\rm der}\succeq aI_T.
 \label{eq:ordinary-penalty-curvature}
\end{equation}
Writing $\overline H_{\rm der}$ for its zero extension, the equivalent
operator formulas are
\begin{equation}
 \boxed{\displaystyle
 A_{\rm der}
 =\bigl[\Pi(\overline H_{\rm der}+X^TDX)\Pi\bigr]^\dagger
 =U\bigl[U^T(\overline H_{\rm der}+X^TDX)U\bigr]^{-1}U^T,}
 \label{eq:ordinary-penalty-inverse}
\end{equation}
where $U$ has orthonormal columns spanning $T$, with value zero if
$T=\{0\}$. Equivalently, $A_{\rm der}$ is the unique symmetric factor
with range $T$ satisfying $A_{\rm der}(P+MJ)=J$.
\end{proposition}

\begin{proof}
The map $G_R$ is globally Lipschitz. The ordinary-derivative defect
test and Borel zero-set selection described in
\cref{prop:ordinary-selection} apply to
$z\mapsto G_R(q+z)$, with candidates $(P,J)$ satisfying
$P=P^T=P^2$ and $J=JP$. Their zero set encodes precisely
\eqref{eq:ordinary-penalty-expansion}. Maximizing the bounded integer
score $\rank P$ selects a largest admissible space measurably.
The dependence of the centered response map on $q$ is continuous
for uniform convergence on compact sets.

Apply subgradient strong monotonicity and cyclic monotonicity to
$G_R(q+Pz)$ and $q+Pz$, and differentiate at zero. This gives the
displayed energy and symmetry identities. The linear-response
curvature argument of \cref{thm:data-response} then gives
\eqref{eq:ordinary-penalty-curvature}. Finally,
$J^T\overline H_{\rm der}J=J^TP=PJ$, so the algebra of
\cref{prop:data-factor-pseudoinverse} proves
\eqref{eq:ordinary-penalty-inverse} and the factor characterization.
\end{proof}

\paragraph{Relationship with the subsequential construction.}

The two formulas \eqref{eq:ordinary-penalty-factor} and
\eqref{eq:ordinary-penalty-inverse} define the same operator.
They do not assert equality with $A_{\rm data}$ from
\cref{def:data-operator}. Ordinary differentiability on $F$ forces
every subsequential response limit to agree with $J_F$ on $F$.
The converse need not hold: a particular subsequence may have a
linear limit on a space where the full derivative fails to exist.

For example, in one dimension take, near zero,
\[
 G(0)=0,\qquad
 G(u)=u\left[2+\tfrac12\sin\!\left(\log\log\frac e{|u|}\right)\right]
 \quad(u\ne0).
\]
On a sufficiently small neighborhood, $G$ is increasing with
derivative bounded above and bounded away from zero. Extend it to
a globally increasing bi-Lipschitz map onto $\mathbb R$ and set
$R(b)=\int_0^b G^{-1}(v)\,dv$. Then $R$ is finite valued and
strongly convex, and $G_R=G$. At $q=0$ the ordinary derivative
does not exist, so $A_{\rm der}=0$. Yet every subsequential
locally uniform limit of $mG(z/m)$ is $z\mapsto cz$ for some
$c\in[3/2,5/2]$. For $M=0$, the original construction therefore
gives $A_{\rm data}=c>0$.

For the prescribed penalties, the situation is simpler.
If $R$ is $C^2$ and strongly convex, $G_R$ is differentiable on the
full space, so
\[
 A_{\rm der}=A_{\rm data}
 =\bigl[\nabla^2R(\widehat b)+X^TDX\bigr]^{-1}=Q_{\rm sm}.
\]
For the positive-weight group Lasso penalty in
\cref{thm:data-group-recovery}, the response has a unique directional
derivative: the quadratic critical-cone variational inequality in
that theorem determines every subsequential limit uniquely.
Lipschitz equicontinuity upgrades this to convergence uniform on
compact sets. A restriction of $G_R$ is ordinarily differentiable
exactly when this directional response is linear on the restriction
space. Thus its maximal differentiability spaces are precisely the
maximal linearity spaces used there, and the same proof yields
\begin{equation}
 \boxed{A_{\rm der}=A_{\rm data}=Q_{\rm grp}.}
 \label{eq:ordinary-penalty-group-recovery}
\end{equation}
This includes inactive boundary groups and overlapping group maps;
no empirical forcing-rank condition is needed.

\paragraph{Statistical guarantees.}

\begin{theorem}[Ordinary penalty-response replacement]
\label{thm:ordinary-penalty-statistics}
Under the hypotheses of \cref{thm:linear} or \cref{thm:single}, all
four respective conclusions hold with $A_{\rm der}$ in place of
$A_{\rm data}$ and the same constants. The pointwise covariance and
rational-weight bounds of \cref{thm:data-weights} also hold.
Almost surely on $\{s\ne0\}$,
\[
 \rank(A_{\rm der}-A_{\rm si})\le4.
\]
\end{theorem}

\begin{proof}
The key observation is that an ordinary differentiable subgradient
graph supplies an admissible penalty-response space. Suppose
$v(\eta)\in\partial R(b(\eta))$, $v(0)=q$, and $b,v$ have derivatives
$B,C$ at zero. Let $F_0=\ran C$ and choose a linear right inverse
$K:F_0\to\mathbb R^d$ of $C$. Since
$b(\eta)=G_R(v(\eta))$ and $G_R$ is Lipschitz,
\[
 G_R(q+z)=G_R(v(Kz))+o(\norm z)
 =b(0)+BKz+o(\norm z),\qquad z\in F_0.
\]
Hence $F_0$ is admissible in
\cref{def:ordinary-penalty-response}.
For the signal-slice jet, $\ran C$ has codimension at most two
almost surely on $\{s\ne0\}$. Maximal dimension therefore gives
$\operatorname{codim}F\le2$ for the selected space.

The selected maps $z\mapsto G_R(q+Pz)$ and $z\mapsto q+Pz$
are themselves an ordinary subgradient graph with derivatives $J,P$.
Comparing this graph with the signal-slice jet on the intersection
of their dual derivative ranges gives the sum-of-codimensions rank
bound, hence rank four, exactly as in the proof of
\cref{thm:ordinary-statistics}. Strong monotonicity in the covariance
metric gives the pointwise energy and denominator bounds.
The existing finite-rank ALO-weight and trace comparisons, followed
by the unconditional prediction and test transfer arguments, give
all four conclusions with the same constants.
\end{proof}

\paragraph{Lean formalization.}
The ordinary penalty-response rule and its statistical endpoints are
checked end to end. \lean{OrdinaryPenaltyResponse.lean} defines
\lean{ordinaryPenaltySelection} and \lean{ordinaryPenaltyALOOperator}
from $R,q,M$, with no convexity modulus, covariance or signal parameter
as a definition input. \lean{exists_ordinaryPenaltySelection} and
\lean{ordinaryPenaltySelection_spec} prove Borel selection, ordinary
Fr\'echet differentiability and maximal dimension.
\lean{hasFDerivAt_penaltyResponse_of_graph_zero} and
\lean{exists_ordinaryPenaltyResponse_of_jet} prove the full ordinary
expansion on a graph jet's dual derivative range;
\lean{ordinaryPenaltySelection_codimension_of_jet} supplies the selected
codimension bound.

\lean{ordinaryPenaltySelection_curvature} proves existence and uniqueness
of the symmetric coercive curvature.
\lean{ordinaryPenaltyALOOperator_isSymmetricFactor} checks the forcing
characterization. The two basis-free formulas
\eqref{eq:ordinary-penalty-factor} and
\eqref{eq:ordinary-penalty-inverse} are checked by
\lean{ordinaryPenaltyALOOperator_eq_responseGram} and
\lean{ordinaryPenaltyALOOperator_eq_symmetricPseudoInverse}.
The former uses the generic singular-Gram identity
\lean{symmetricFactor_eq_responseGram}, including zero responses.
\lean{ordinaryPenaltyALOOperator_eq_smooth_inverse} and
\lean{ordinaryPenaltyALOOperator_eq_groupLassoEmpiricalALOOperator}
prove exact recovery of the prescribed matrices.
\lean{ordinaryPenaltyALOOperator_eq_dataDriven_of_groupLasso} also checks
agreement with the original Section~\ref{sec:data-operator} rule,
including overlapping and inactive boundary groups.

\lean{ordinaryPenaltyALOOperator_control},
\lean{norm_whiten_ordinaryPenaltyALOOperator_le} and
\lean{ordinaryPenaltyALOWeight_bound} give the pointwise bounds.
\lean{ae_covarianceOrdinaryPenaltySignalOperator_rank} proves the
unconditional rank comparison. The four final model endpoints are
\lean{regression_ordinaryPenalty_ALO_unconditional},
\lean{regression_ordinaryPenalty_meanField_unconditional},
\lean{singleIndex_ordinaryPenalty_ALO_unconditional}, and
\lean{singleIndex_ordinaryPenalty_meanField_unconditional}.
They include integrability and expected absolute test-error bounds,
with the same constants as the original endpoints and no extra
selector, differentiability or codimension premises.
\lean{singleIndex_observedOrdinaryPenaltyOperator_eq_original_loss}
identifies their operator with the selection from the original penalty,
observed scores and loss curvature. No new axiom or proof placeholder
is used. The illustrative oscillatory-response counterexample above
is a mathematical example and is not a separate Lean construction.

\section{An operator selected from ordinary empirical derivatives}
\label{sec:empirical-derivative}

There is another data-driven construction, closer to
\cref{def:operator}, which uses ordinary derivatives of the empirical
fit rather than a selected subsequential penalty response. It need
not equal \cref{def:data-operator}. All construction inputs below are
in the original coordinates; whitening is used only for analysis.

Hold the observed responses fixed and, for a variable design $W$, write
\begin{equation}
 b_y(W)=\argmin_b\left\{\sum_{i=1}^n L_{y_i}(w_i^Tb)+R(b)\right\},\qquad
 q_y(W)=-\sum_{i=1}^n L_{y_i}'(w_i^Tb_y(W))w_i.
 \label{eq:ordinary-empirical-maps}
\end{equation}
At the observed $X$, retain $s,D,M$ from \eqref{eq:data-inputs} and set
\begin{equation}
 \mathcal L_X(E)=-E^Ts-X^TDE\widehat b,
 \qquad E\in\mathbb R^{n\times p}.
 \label{eq:ordinary-forcing}
\end{equation}

\begin{definition}[Selected empirical design derivatives]
\label{def:ordinary-design}
A linear design space $\mathcal E\subseteq\mathbb R^{n\times p}$ is
\emph{admissible at $(X,y)$} if the restrictions of $b_y,q_y$ to
$X+\mathcal E$ are Fr\'echet differentiable at $X$, with derivatives
$B_{\mathcal E},C_{\mathcal E}$, and
\begin{equation}
 C_{\mathcal E}=\mathcal L_X|_{\mathcal E}-MB_{\mathcal E}.
 \label{eq:ordinary-kkt}
\end{equation}
Among admissible spaces, select one maximizing
$\rank(\mathcal L_X|_{\mathcal E})$. Among these, maximize
$\dim\mathcal E$ and resolve remaining ties by a fixed Borel rule.
Write $B,C$ for its derivatives and $T=\ran B$.
Define $A_{\rm emp}$ by
\begin{equation}
 \boxed{\quad A_{\rm emp}^T=A_{\rm emp},\qquad
 A_{\rm emp}\mathcal L_X|_{\mathcal E}=B,\qquad
 \ran A_{\rm emp}=T.\quad}
 \label{eq:ordinary-factor}
\end{equation}
The zero space is admissible. The selections use only the observed
losses, the penalty and $X$, with no covariance or signal parameter.
\end{definition}

\begin{proposition}[Selection and original-coordinate restricted inverse]
\label{prop:ordinary-selection}
The selections in \cref{def:ordinary-design} can be made Borel, and
\eqref{eq:ordinary-factor} has exactly one solution. There is a unique
symmetric operator $H:T\to T$ with
\begin{equation}
 HB(E)=P_TC(E),\qquad H\succeq aI_T.
 \label{eq:ordinary-curvature}
\end{equation}
If $U_T$ has orthonormal columns spanning $T$ and $H_T$ is the matrix
of $H$ in that basis, then
\begin{equation}
 \boxed{\displaystyle A_{\rm emp}
 =U_T\bigl[H_T+U_T^TX^TDXU_T\bigr]^{-1}U_T^T
 =\bigl[P_T(\overline H+X^TDX)P_T\bigr]^\dagger.}
 \label{eq:ordinary-restricted-inverse}
\end{equation}
Here $\overline H$ is the zero extension and $A_{\rm emp}=0$ for
$T=\{0\}$. If the full design space is admissible, it is selected.
\end{proposition}

\begin{proof}
For fixed observed losses, $b_y$ and $q_y$ are locally Lipschitz in
the design. Encode a candidate by its orthogonal projector $P$ and
the ambient extensions $BP,CP$ of its two derivatives. Restricted
differentiability is equivalently the ordinary uniform expansion
\[
 b_y(X+tPE)=b_y(X)+tBE+o(|t|),\qquad
 q_y(X+tPE)=q_y(X)+tCE+o(|t|),
\]
uniformly over $\norm E\le1$, together with $B=BP,C=CP$.
Here is a Borel selection argument using only closed relations.
For a centered continuous empirical map $f$, put
\[
 \begin{aligned}
 e_m(f;P,B)&=\sup_{\norm E\le1}
   \norm{2^m\{f(2^{-m}PE)-f(0)\}-BE},\\
 \Delta(f;P,B)&=\limsup_{m\to\infty}\min\{1,e_m(f;P,B)\}.
 \end{aligned}
\]
The defect is zero exactly when the ordinary derivative of $f\circ P$
at zero is $B$. Dyadic scales suffice because every sufficiently small
increment lies between two consecutive such scales. The defect is
Borel in the empirical map and candidate matrices, and is continuous
in the candidate matrices for each fixed locally Lipschitz map.
For the last assertion one may extend the map from a neighborhood of
zero to a globally Lipschitz finite-dimensional map; the defect is
unchanged, and its Lipschitz bound in $(P,B)$ is then uniform in scale.
This is a full-sequence ordinary-derivative test, with no selected limit.

Add the two defects and the nonnegative norm errors for the projector,
the identities $B=BP,C=CP$, and \eqref{eq:ordinary-kkt}. The resulting
continuous-in-candidate Borel error has precisely the admissible
triples $(P,B,C)$ as its zero set, and the zero triple is always a zero.
Select a zero in each closed ball with center in a fixed countable
dense set and radius $1/(k+1)$, using a fixed default if the ball has
no zero. Each such selection is Borel: its graph is a closed relation
with compact fibers, and fixed dense-sequence distances provide unique
lexicographic tie-breaking. These selections form a countable dense
family in every zero fiber.

On the finite-dimensional candidate space, rank upper bounds are
closed. Consequently the bounded integer score
\[
 (np+1)\rank(\mathcal L_XP)+\rank P
\]
has open upper level sets and a dense family hits a maximizer whenever
one exists. Select the first member attaining the largest score.
Countable tests show this choice is Borel. Since $\rank P\le np$,
the score maximizes forcing rank first and design dimension second,
as required. The observed score functions are themselves Borel
parameters in the topology of uniform convergence on compact sets.

The subgradient graph relation $q_y(W)\in\partial R(b_y(W))$
implies $\langle B(E),C(E)\rangle\ge a\norm{B(E)}^2$ and the
mixed symmetry of the derivative relation. The same tangent-curvature
argument as in \cref{thm:operator-existence} gives
\eqref{eq:ordinary-curvature}. Projecting \eqref{eq:ordinary-kkt}
onto $T$ gives \eqref{eq:ordinary-restricted-inverse}, including
the factor equation and its uniqueness. Finally, restriction cannot
increase forcing rank, and the full space has the largest dimension;
thus the tie-breaking rule selects it whenever it is admissible.
\end{proof}

The form \eqref{eq:ordinary-restricted-inverse} is exactly the
original-coordinate counterpart of \eqref{eq:restricted-inverse}.
It has the same algebraic form as \eqref{eq:data-restricted-inverse},
but its $T,H$ are supplied by empirical design derivatives. They need
not equal $T_q,H_q$. For instance, with
$R(b)=a\norm b^2/2$ and $L_{y_i}(t)=\log\cosh t$, the fit is zero
for every design. Then $A_{\rm emp}=0$, whereas
$A_{\rm data}=(aI+X^TX)^{-1}$; all fitted scores are zero.

\begin{example}[Why the full design space need not be admissible]
\label{ex:ordinary-threshold}
Take $n=p=1$, $R(b)=b^2/2+|b|$, $L_y(t)=-yt$, and
\[
 y=f(x)=\begin{cases}1/x,&x\ge1,\\0,&x<1.\end{cases}
\]
This is a measurable single-index response with $\beta=1$; its
observed loss score is constant, bounded by one, and 1-Lipschitz.
Holding the realized $y$ fixed, the empirical fit at a variable design
$w$ is
\[
 b_y(w)=\operatorname{soft}(wy,1),\qquad
 \operatorname{soft}(u,1)=\sign(u)(|u|-1)_+.
\]
On $\{x\ge1\}$, $xy=1$ and $y>0$, so the left derivative of
$b_y$ at $w=x$ is zero and the right derivative is $y$. Thus the
full design space is not admissible on an event of positive Gaussian
probability, even though for each fixed $y$ the fit is differentiable
at almost every $w$. The realized response can put the design exactly
at a kink. Setting the full-space derivative to zero there does not
automatically give the rank comparison required for statistical transfer.
\end{example}

\begin{theorem}[Comparison and statistical replacement]
\label{thm:ordinary-statistics}
Under the hypotheses of \cref{thm:linear} or \cref{thm:single},
the operator $A_{\rm emp}$ is data-driven and satisfies all four
respective conclusions with $A_{\rm data}$ replaced by $A_{\rm emp}$
and the same constants. In particular these are unconditional ALO
and mean-field bounds. Almost surely on $\{s\ne0\}$,
\begin{equation}
 \rank(A_{\rm emp}-A_{\rm si})\le4,
 \label{eq:ordinary-rank}
\end{equation}
where $A_{\rm si}$ is the canonical signal-slice operator.
The bounds of \cref{thm:data-weights} also hold for $A_{\rm emp}$.
\end{theorem}

\begin{proof}
On a fixed signal slice, the empirical fit and its subgradient have
ordinary derivatives almost everywhere, and the Lipschitz-score chain
rule gives \eqref{eq:ordinary-kkt}. Thus the true signal-preserving
space $\{E:E\beta=0\}$ is an admissible competitor almost surely.
This is also valid for $\beta=0$, when the competitor is the full space.
If $s\ne0$, its forcing range has codimension at most two. The
selected forcing range therefore has codimension at most two as well.
For either jet, strong monotonicity and the factor equation give
\[
 \ker C=\ker(\mathcal L_X|_{\mathcal E}),\qquad
 \rank C=\rank(\mathcal L_X|_{\mathcal E}).
\]
Indeed $C(E)=0$ forces $B(E)=0$ by strong monotonicity and then
$\mathcal L_X(E)=0$; conversely, $\mathcal L_X(E)=0$ forces
$B(E)=0$ by the factor equation and then $C(E)=0$ by KKT.

For two differentiable primal--dual graphs at the same dual point,
the coefficient derivatives agree whenever the subgradient derivatives
agree. One can prove this directly by replacing the nonlinear dual
increment by its derivative and using the Lipschitz inverse
subgradient $G_R$. Thus on the intersection of the two dual derivative
ranges, of codimension at most four, both responses have the same
value $Kz$. Their factor equations give
\[
 A_{\rm emp}(z+MKz)=Kz=A_{\rm si}(z+MKz).
\]
The map $z\mapsto z+MKz$ is injective there: if its value is zero,
then $a\norm{Kz}^2\le\langle Kz,z\rangle
=-\langle Kz,MKz\rangle\le0$, so $Kz=z=0$.
Rank--nullity proves \eqref{eq:ordinary-rank}.

The restricted inverse supplies positivity and the energy bounds.
Using $R\circ S^{-1}$ for the proof also gives
$\norm{SA_{\rm emp}S^T}_{\rm op}\le(n\mu)^{-1}$ and the same
quantitative denominator bounds as for $A_{\rm data}$.
The Gaussian finite-rank rational-weight comparison therefore has
the same rank-four constant; the trace comparison has bound
$64/(n\mu^2)$. Applying the existing unconditional prediction and
test transfer arguments proves the four conclusions. When $s=0$,
the score-weighted comparisons vanish. All conditioning used here is
on the noise and signal coordinates, followed by integration; no fitted
quantity is conditioned on.
\end{proof}

\paragraph{Lean formalization.}
The construction and the four statistical conclusions of
\cref{thm:ordinary-statistics} are checked end to end.
\lean{ProjectionDerivativeDefect.lean} proves the dyadic ordinary-derivative
test and its continuity in the candidate matrices.
\lean{ContinuousZeroSelection.lean} constructs a countable Borel family
dense in each zero fiber, then selects a maximizer of a bounded integer
score with open upper level sets.
\lean{exists_measurable_maximal_empiricalDesign} applies these results to
both rank maximizations. The rule is fixed on the canonical observed-data
space by \lean{exists_ordinaryEmpiricalSelection};
\lean{ordinaryEmpiricalOperator} uses only the observed score functions,
design and penalty. Its measurability, symmetric-factor identity and
full-space selection when admissible are checked in
\lean{OrdinaryEmpiricalOperator.lean}.

\lean{EmpiricalAffineCompetitor.lean} transfers ordinary derivatives on a
signal slice to an admissible original-coordinate design space.
\lean{ae_covarianceOrdinaryEmpiricalSignalOperator_rank} proves the resulting
almost-sure rank-four comparison when the signal codimension is at most
one; the norm and denominator bounds follow from
\lean{HasEmpiricalResponseJet.rational_control}.
The final theorems are
\lean{regression_ordinaryEmpirical_ALO_unconditional},
\lean{regression_ordinaryEmpirical_meanField_unconditional},
\lean{singleIndex_ordinaryEmpirical_ALO_unconditional}, and
\lean{singleIndex_ordinaryEmpirical_meanField_unconditional}.
They include integrability and expected absolute test-error bounds under
the original covariance-design law and arbitrary independent noise,
with the same constants as the penalty-response endpoints. No
measurable-selector, derivative-existence or rank-comparison premise is
added to the model assumptions.
\lean{singleIndex_observedOrdinaryEmpiricalOperator_eq} identifies the
operator in these bounds with the original-coordinate observed-data rule.
No new axiom or proof placeholder is used. The explicit matrix
pseudoinverse display in \eqref{eq:ordinary-restricted-inverse}, as distinct
from the symmetric-factor characterization, is proved algebraically above.

\subsection{When the full design space is admissible almost surely}
\label{sec:ordinary-full-space}

The full design space can be used almost surely under a condition on
the joint data law, without any additional smoothness of the penalty
or loss scores. Write $Y=(y_1,\ldots,y_n)$ and let
$\lambda_{np}$ denote Lebesgue measure on $\mathbb R^{n\times p}$.

\begin{proposition}[Conditional absolute continuity]
\label{prop:ordinary-full-space}
Under the standing assumptions of this section, suppose that $Y$ takes
values in a standard Borel space and
\begin{equation}
 \mathcal L(X\mid Y=y)\ll\lambda_{np}
 \qquad\text{for }\mathcal L(Y)\text{-almost every }y.
 \label{eq:ordinary-conditional-density}
\end{equation}
Then $\mathcal E=\mathbb R^{n\times p}$ is admissible at $(X,Y)$
almost surely. Consequently the selection in
\cref{def:ordinary-design} chooses the full space almost surely, and
its derivatives are the ordinary full-design derivatives
$B=Db_Y(X)$ and $C=Dq_Y(X)$.
\end{proposition}

\begin{proof}
Fix $y$. Both maps in \eqref{eq:ordinary-empirical-maps} are locally
Lipschitz. Rademacher's theorem therefore gives their Fr\'echet
derivatives outside a Lebesgue-null set
of designs. Write $s_y(W),D_y(W)$ for the scores and curvatures in
\eqref{eq:data-inputs}, evaluated at $(W,b_y(W))$.
The almost-everywhere chain rule for Lipschitz scores,
with the zero convention for their derivatives at kinks, gives
\[
 Dq_y(W)[E]
 =-E^Ts_y(W)-W^TD_y(W)E b_y(W)
   -W^TD_y(W)W\,Db_y(W)[E]
\]
for Lebesgue-almost every $W$. Thus the full space satisfies both
differentiability and \eqref{eq:ordinary-kkt} outside a null set $N_y$.
The admissibility test used in \cref{prop:ordinary-selection} makes
the joint exceptional set Borel. Conditional absolute continuity gives
$\mathbb P(X\in N_y\mid Y=y)=0$ for almost every $y$; integrating
proves full-space admissibility almost surely. The final assertion
follows from \cref{prop:ordinary-selection}.
\end{proof}

In particular, \eqref{eq:ordinary-conditional-density} holds whenever
$X$ has a density and $Y$ has countable range, including binary or
finite-label response vectors. Indeed, for every $y$ with
$\mathbb P(Y=y)>0$ and every Lebesgue-null set $N$,
\[
 \mathbb P(X\in N\mid Y=y)
 =\frac{\mathbb P(X\in N,\,Y=y)}{\mathbb P(Y=y)}=0.
\]
No independence of $X$ and $Y$ is required.

\begin{example}[Linear regression with independent noise having a density]
\label{ex:ordinary-linear-density}
Consider $Y=X\beta+\varepsilon$, where $X$ has a density $f_X$ and
the noise vector $\varepsilon$ is independent of $X$ with a Lebesgue
density $f_\varepsilon$ on $\mathbb R^n$. The joint density is
\[
 f_{X,Y}(x,y)=f_X(x)f_\varepsilon(y-x\beta).
\]
Consequently, for $\mathcal L(Y)$-almost every $y$,
\[
 f_{X\mid Y=y}(x)
 =\frac{f_X(x)f_\varepsilon(y-x\beta)}{f_Y(y)},\qquad
 f_Y(y)=\int f_X(w)f_\varepsilon(y-w\beta)\,dw.
\]
Thus \eqref{eq:ordinary-conditional-density} holds, and
\cref{prop:ordinary-full-space} permits the full design space almost
surely. In particular, this applies to the Gaussian design in
\cref{thm:linear} with independent Gaussian, Laplace or uniform noise
coordinates. Neither smoothness nor strict positivity of the noise
density is required; dependence among noise coordinates is allowed
provided the vector has a joint density. No additional smoothness of
$R$ or of the loss scores is imposed, and the ALO and mean-field
conclusions of \cref{thm:ordinary-statistics} hold with the same constants.
\end{example}

In contrast, in \cref{ex:ordinary-threshold}, conditioning on a
positive response $Y=y$ fixes $X=1/y$, so the conditional design law
is singular and can charge the kink. Marginal Gaussianity of $X$
alone does not exclude this behavior. Under
\eqref{eq:ordinary-conditional-density}, the construction and all
statistical conclusions of \cref{thm:ordinary-statistics} therefore
use the full-space derivatives almost surely; the existing Borel
selection specifies the operator on the exceptional null set as well.

\section{Unconditional data-driven ALO and mean-field bounds}
\label{sec:unconditional}

Let $n\ge1$, $\mu,\delta>0$, and assume $p\le n/\delta$.
The rows of $X$ are independent $N(0,\Sigma)$ vectors with
$\Sigma=S^TS\succ0$. Throughout this section, $R$ is finite valued,
strongly convex, minimized at zero, and
$\widetilde R(\theta)=R(S^{-1}\theta)$ is $n\mu$-strongly convex.
Equivalently, the last condition requires convexity of
$R(b)-n\mu b^T\Sigma b/2$.

For a loss $L_y(t)$, define the full and deleted fits using the
\emph{same} penalty $R$:
\begin{equation}
 \widehat b=\argmin_b\left\{\sum_{j=1}^n L_{y_j}(x_j^Tb)+R(b)\right\},
 \qquad
 \widehat b^{(-i)}
 =\argmin_b\left\{\sum_{\substack{j=1\\j\ne i}}^nL_{y_j}(x_j^Tb)+R(b)\right\}.
 \label{eq:full-loo}
\end{equation}
Write $\psi_y=L_y'$, and use the zero convention for
$d_i=\psi_{y_i}'(x_i^T\widehat b)$ at score kinks.
In both theorems below, compute $A_{\rm data}$ from
\cref{def:data-operator} using the original penalty $R$, the observed
subgradient $q=-\sum_{i=1}^n\psi_{y_i}(x_i^T\widehat b)x_i$, and $M=X^TDX$.
The same operator is used in ALO and in the mean-field prediction.
For a matrix $A$, write
\begin{align}
 t_i^{\rm ALO}(A)
   &=x_i^T\widehat b+
     \psi_{y_i}(x_i^T\widehat b)
       \frac{x_i^TAx_i}{1-d_i x_i^TAx_i},
       \label{eq:alo-prediction}\\
 t_i^{\rm MF}(A)
   &=x_i^T\widehat b+
     \psi_{y_i}(x_i^T\widehat b)\tr(A\Sigma),&
 t_i^{\rm LOO}&=x_i^T\widehat b^{(-i)},\label{eq:mf-prediction}\\
 \alo(A)&=\frac1n\sum_{i=1}^n g(t_i^{\rm ALO}(A),y_i),&
 \loo&=\frac1n\sum_{i=1}^n g(t_i^{\rm LOO},y_i),\label{eq:test-criteria}\\
 \operatorname{MF}(A)&=\frac1n\sum_{i=1}^n g(t_i^{\rm MF}(A),y_i).
 \label{eq:mf-criterion}
\end{align}
For $A=A_{\rm data}$ every denominator in \eqref{eq:alo-prediction}
is strictly positive by \cref{thm:data-weights}.
The test function $g:\R\times\mathcal Y\to\R$ is jointly measurable and
satisfies
\begin{equation}
 |g(t,y)-g(t',y)|\le |t-t'|(1+|t|+|t'|)
 \quad\text{for every }t,t',y.
 \label{eq:test-condition}
\end{equation}

Let $C_{\rm ALO}$, $C_{\rm emp}$ and $H_{\rm pred}$ be the finite
nonnegative functions of $(\mu,\delta)$ specified below. Set
\begin{align}
 C_{\rm MF}(n;\mu,\delta)
   &=2C_{\rm emp}(\mu,\delta)+\frac{128}{n\mu^2},
     \label{eq:data-mf-constant}\\
 \mathcal T(C;\mu,\delta)
   &=\sqrt{C}\sqrt{3+18\bigl(H_{\rm pred}(\mu,\delta)
                      +C_{\rm emp}(\mu,\delta)\bigr)+6C}.
     \label{eq:data-test-constant}
\end{align}
For fixed $\mu,\delta$, $C_{\rm MF}$ is bounded uniformly over $n\ge1$.
Thus each test bound below is of order $n^{-1/2}$.

\begin{theorem}[Data-driven linear regression]
\label{thm:linear}
Assume the preceding design, penalty, and test-function conditions.
Let
$y_i=x_i^T\beta+\eps_i$, where $\beta\in\R^p$ is fixed and the
entire vector $\eps$ is independent of $X$, with an arbitrary
probability law. Let $L_y(t)=\rho(y-t)$, where $\rho$ is convex
and differentiable, $\rho'$ is 1-Lipschitz, and
$|\rho'(t)|\le1$ for every $t$.

All four nonnegative random quantities below are integrable, and
\begin{align}
 \E\sum_{i=1}^n
   |t_i^{\rm ALO}(A_{\rm data})-t_i^{\rm LOO}|^2
       &\le C_{\rm ALO}(\mu,\delta),\label{eq:linear-squared}\\
 \E\bigl[|\alo(A_{\rm data})-\loo|\bigr]
       &\le\frac{\mathcal T(C_{\rm ALO}(\mu,\delta);\mu,\delta)}{\sqrt n},
       \label{eq:linear-unconditional}\\
 \E\sum_{i=1}^n
   |t_i^{\rm MF}(A_{\rm data})-t_i^{\rm LOO}|^2
       &\le C_{\rm MF}(n;\mu,\delta),\label{eq:linear-mf-squared}\\
 \E\bigl[|\operatorname{MF}(A_{\rm data})-\loo|\bigr]
       &\le\frac{\mathcal T(C_{\rm MF}(n;\mu,\delta);\mu,\delta)}{\sqrt n}.
       \label{eq:linear-mf-test}
\end{align}
The operator $A_{\rm data}$ and the ALO criterion are determined by
$(y,X)$, the loss and the penalty; neither $\beta$, $\eps$ nor
$\Sigma$ is an input to their construction. The mean-field criterion
additionally uses $\Sigma$ through $\tr(A_{\rm data}\Sigma)$.
\end{theorem}

\paragraph{Lean formalization.}
The two endpoints are
\lean{regression_dataDriven_ALO_unconditional} and
\lean{regression_dataDriven_meanField_unconditional}.
They allow arbitrary positive-definite covariance and an arbitrary
joint noise law independent of the design. No $C^2$ assumption is
made on the loss or penalty.

\begin{theorem}[Data-driven single-index models]
\label{thm:single}
Assume the design, penalty and test-function conditions at the start
of this section. Let $U$ be an arbitrary random element independent
of $X$, and let
\begin{equation}
 y_i=f_i(U,x_i^T\beta),\qquad i=1,\ldots,n,
 \label{eq:single-index}
\end{equation}
where $\beta\in\R^p$ is fixed and each $f_i$ is measurable.
For every $y$, assume that $L_y$ is convex and differentiable,
$\psi_y=L_y'$ is 1-Lipschitz, and $|\psi_y(t)|\le1$ for every $t$.
Assume also that $y\mapsto\psi_y(t)$ is measurable for each $t$.

With $A_{\rm data}$ constructed from the realized data as above,
all four quantities below are integrable, and
\begin{align}
 \E\sum_{i=1}^n
   |t_i^{\rm ALO}(A_{\rm data})-t_i^{\rm LOO}|^2
       &\le C_{\rm ALO}(\mu,\delta),\label{eq:single-squared}\\
 \E\bigl[|\alo(A_{\rm data})-\loo|\bigr]
       &\le\frac{\mathcal T(C_{\rm ALO}(\mu,\delta);\mu,\delta)}{\sqrt n},
       \label{eq:single-unconditional}\\
 \E\sum_{i=1}^n
   |t_i^{\rm MF}(A_{\rm data})-t_i^{\rm LOO}|^2
       &\le C_{\rm MF}(n;\mu,\delta),\label{eq:single-mf-squared}\\
 \E\bigl[|\operatorname{MF}(A_{\rm data})-\loo|\bigr]
       &\le\frac{\mathcal T(C_{\rm MF}(n;\mu,\delta);\mu,\delta)}{\sqrt n}.
       \label{eq:single-mf-test}
\end{align}
The constants do not depend on $\beta$, the links $f_i$, the noise
law, or the covariance beyond the stated penalty and aspect-ratio
conditions. The case $\beta=0$ is included.
The operator and ALO criterion do not require knowledge of $U$,
$\beta$, the links $f_i$, or $\Sigma$; the realized responses and
the loss family suffice. The mean-field criterion uses $\Sigma$.
\end{theorem}

\paragraph{Lean formalization.}
The two endpoints are
\lean{singleIndex_dataDriven_ALO_unconditional} and
\lean{singleIndex_dataDriven_meanField_unconditional}.
Their shared general-covariance proofs are
\lean{dataDrivenCovariance_ALO_unconditional} and
\lean{dataDrivenCovariance_meanField_unconditional}.
The formalization encodes a loss by integrating its bounded Lipschitz
score. The identities
\lean{singleIndex_covarianceEstimator_eq_original_loss},
\lean{singleIndex_covarianceLOO_eq_original_loss}, and
\lean{singleIndex_observedDataDrivenOperator_eq_original_loss}
identify the full fit, deleted fits and $A_{\rm data}$ with those of
the original losses. Additive constants in $L_y$ have no effect.
The regression theorem is the specialization
$f_i(U,t)=t+\eps_i$, with $U=\eps$ and
$\psi_y(t)=-\rho'(y-t)$.

\begin{corollary}[Group penalties with the data-driven operator]
\label{cor:data-group}
For either model, take the ridge-augmented group penalty
\begin{equation}
 R(b)=\sum_{k=1}^K\lambda_k\norm{P_kb}
             +\frac a2\norm b^2,\qquad
 \lambda_k\ge0,\quad a>0,\quad a\ge n\mu\norm{\Sigma}_{\rm op}.
 \label{eq:group-penalty}
\end{equation}
Then all four conclusions of \cref{thm:linear,thm:single}, as applicable,
hold with this penalty's $A_{\rm data}$ and the same constants.
The linear maps $P_k:\R^p\to\R^p$ may overlap and need not jointly
cover the coefficient space.
If all retained weights are positive, $A_{\rm data}=Q_{\rm grp}$ by
\cref{thm:data-group-recovery}, so these same conclusions hold for
the prescribed group matrix, without a coverage condition.
For zero-weight terms, remove those groups before forming $Q_{\rm grp}$.
\end{corollary}

\paragraph{Lean formalization.}
The single-index specializations are
\lean{singleIndex_groupLasso_dataDriven_ALO_unconditional} and
\lean{singleIndex_groupLasso_dataDriven_meanField_unconditional};
regression is included by the response specialization above.
These results use $A_{\rm data}$, which agrees pointwise with
$Q_{\rm grp}$ for positive retained weights.
The older independently proved results for $Q_{\rm grp}$ are stated
in \cref{cor:prescribed-group}.

\paragraph{Constants in the checked bounds.}
One may take $C_{\rm ALO}$, $C_{\rm emp}$ and $H_{\rm pred}$ to be
\lean{dataDrivenSignalALOConstant},
\lean{empiricalMeanFieldConstant} and
\lean{meanFieldPredictionMomentConstant}, respectively, evaluated
at $(\mu,\delta)$. More explicitly,
\begin{equation}
 C_{\rm ALO}
 =3C_{\rm rank}+3C_{\rm weight}+3C_{\rm emp}
       +\frac{4C_{\rm bad}}{\mu^2},
 \label{eq:data-alo-constant}
\end{equation}
Here $C_{\rm rank}$ is \lean{gaussianRankWeightConstant}
evaluated at $(4,\mu,\delta)$; $C_{\rm weight}$ is
\lean{constructedSignalWeightConstant} at $(\mu,\delta)$; and
$C_{\rm bad}$ is \lean{gaussianBadRowFourthConstant} at $\delta$.
The quantity $\mathcal T(C;\mu,\delta)/\sqrt n$ is exactly
\lean{dataDrivenTestRate} evaluated at $(n,C,\mu,\delta)$.
In particular the covariance extension introduces no additional
condition-number factor in these constants.

\begin{remark}[Meaning of the expectations]
All expectations in \cref{thm:linear,thm:single} are under the full
original noise/design law. For the test criteria the absolute value
is inside the expectation and outside the average over observations.
There is no indicator $I_\Omega$, division by an event probability, or
conditioning on the estimator. Noise coordinates may be dependent
and no noise moments are assumed. It is the \emph{difference} of the
two criteria that is proved integrable; the criteria separately need
not be integrable.
\end{remark}

\paragraph{Proof of the unconditional bounds.}
The comparison with the derivative-based operator is made on fixed
signal slices, then integrated under the original law. As detailed
in \cref{thm:data-comparison}, the two whitened operators differ by
rank at most four wherever the fitted score vector is nonzero.
Their weighted rational-weight discrepancy has unconditional mean
square at most $C_{\rm rank}$. Their summed squared, score-weighted
trace discrepancy is at most $64/(n\mu^2)$, giving \eqref{eq:data-mf-constant} after a
squared triangle inequality with the old mean-field error.

For ALO, the old weight-to-trace estimate is used on
$\Omega=\{\max_i\norm{z_i}\le\sqrt p+\sqrt n\}$.
The pointwise bound \eqref{eq:data-weight-bound} and deletion
stability give, on every sample,
\begin{equation}
 |t_i^{\rm ALO}(A_{\rm data})-t_i^{\rm LOO}|
       \le\frac{2\norm{z_i}^2}{n\mu}.
 \label{eq:global-envelope}
\end{equation}
Gaussian tails imply
$\E[1_{\Omega^c}\sum_{i=1}^n\norm{z_i}^4]\le n^2C_{\rm bad}(\delta)$.
A three-term squared bound on $\Omega$, together with
\eqref{eq:global-envelope} on $\Omega^c$, yields
\eqref{eq:data-alo-constant}. The old mean-field prediction energy
bound and the test condition \eqref{eq:test-condition} give
\eqref{eq:data-test-constant}. Finally, covariance transport and
integration over arbitrary independent noise preserve all bounds.
The shared Lean arguments are
\lean{signalOperator_ALO_error_unconditional},
\lean{signalOperator_ALO_test_unconditional}, and
\lean{signalOperator_meanField_unconditional}.

\section{Relationship with prescribed inverse matrices}
\label{sec:comparison}

The operator $A_{\rm data}$ in the statistical theorems is defined
by \cref{def:data-operator}. For its proof, let $A_{\rm si}$ denote
the canonical operator of \cref{def:operator} on the single-index
slice: take $c=0$, hold $U$ and $Z\zeta=X\beta$ fixed with
$\zeta=S\beta$, and differentiate along
$\mathcal E=\{E:E\zeta=0\}$ with $\ell_i=L_{y_i}$.
Regression is included by taking $U=\eps$ and $f_i(U,t)=t+\eps_i$.
The comparison operator $A_{\rm si}$ is used only in the proof.
The usual inverse Hessian and active-group inverse instead use a
prescribed coefficient space. We state the relationships separately.

\begin{definition}[Prescribed matrices]
\label{def:group-inverse}
For a twice continuously differentiable penalty, put
\begin{equation}
 Q_{\rm sm}=\bigl(\nabla^2R(\widehat b)+X^TDX\bigr)^{-1}.
 \label{eq:smooth-inverse}
\end{equation}
For \eqref{eq:group-penalty}, define
\begin{align}
 T_{\rm grp}
   &=\{z:P_kz=0\text{ whenever }P_k\widehat b=0\},\\
 K(u)&=
   \begin{cases}
    \norm u^{-1}\bigl(I_p-uu^T/\norm u^2\bigr),&u\ne0,\\
    0,&u=0,
   \end{cases}\\
 H_{\rm grp}
   &=aI_p+\sum_{k=1}^K\lambda_k P_k^TK(P_k\widehat b)P_k.
\end{align}
If $V_{\rm grp}$ has orthonormal columns spanning $T_{\rm grp}$, let
\begin{equation}
 Q_{\rm grp}
 =V_{\rm grp}
   \bigl[V_{\rm grp}^T(H_{\rm grp}+X^TDX)V_{\rm grp}\bigr]^{-1}
   V_{\rm grp}^T.
 \label{eq:group-inverse}
\end{equation}
Set $Q_{\rm grp}=0$ when $T_{\rm grp}=\{0\}$.
For disjoint coordinate groups this is the active-block inverse,
extended by zero off the active coordinates. The definition also
allows overlapping linear group maps.
\end{definition}

\begin{theorem}[Exact recovery of the smooth matrix]
\label{thm:data-smooth-recovery}
Let $R:\R^p\to\R$ be twice continuously differentiable and
$a$-strongly convex, with $a>0$. For every $b\in\R^p$,
$q\in\partial R(b)$, and symmetric positive semidefinite $M$,
\[
 A_{\rm data}(R,q,M)=\bigl(\nabla^2R(b)+M\bigr)^{-1}.
\]
In particular, at every empirical fit,
$A_{\rm data}=Q_{\rm sm}$. This identity requires neither a
nonzero fitted score nor surjectivity of a design derivative.
It holds in the original coordinates for every design covariance.
Consequently the ALO and mean-field conclusions in
\cref{sec:unconditional} hold with $Q_{\rm sm}$ whenever their
statistical hypotheses and the stated smoothness hypothesis hold,
with the same constants.
\end{theorem}

\begin{proof}
Use the primal--dual graph $b(\eta)=\eta$,
$q(\eta)=\nabla R(\eta)$. Its derivatives are $I_p$ and
$H=\nabla^2R(b)$, respectively. Strong convexity gives
$H\succeq aI_p$, so the dual derivative is surjective.
The response comparison therefore has codimension zero and rank
zero. The matrix $(H+M)^{-1}$ satisfies its forcing equation,
$(H+M)^{-1}(H\eta+M\eta)=\eta$, proving the identity.
\end{proof}

Lean checks this in
\lean{dataDrivenALOOperator_eq_of_surjective_jet},
\lean{dataDrivenALOOperator_eq_smooth_inverse}, and the empirical
specialization
\lean{dataDrivenALOOperator_eq_smoothEmpiricalALOOperator}.

\begin{theorem}[Exact recovery of the group matrix]
\label{thm:data-group-recovery}
Let
\[
 R(z)=\sum_{k=1}^K\lambda_k\norm{P_kz}+\frac a2\norm z^2,
 \qquad a>0,\quad \lambda_k>0.
\]
For every $b\in\R^p$, $q\in\partial R(b)$, and symmetric positive
semidefinite $M$, define $T=T_{\rm grp}$ and $H=H_{\rm grp}$ from
\cref{def:group-inverse}, with $\widehat b$ replaced by $b$.
Then the selected response satisfies
\[
 \ran J=T,\qquad \Pi_T HJ=\Pi_T P_F,
\]
and
\[
 A_{\rm data}(R,q,M)
   =V_{\rm grp}\bigl[V_{\rm grp}^T(H+M)V_{\rm grp}\bigr]^{-1}
      V_{\rm grp}^T
   =\bigl[\Pi_T(H+M)\Pi_T\bigr]^\dagger.
\]
Thus $A_{\rm data}=Q_{\rm grp}$ at every empirical fit.
This includes empty active sets, overlapping linear group maps,
and inactive groups at the boundary of their KKT conditions.
No strict complementarity, group coverage, nonzero-score condition,
or covariance assumption is needed for the identity.
The ALO weights and $\tr(A_{\rm data}\Sigma)$ coincide exactly with
those of $Q_{\rm grp}$.
\end{theorem}

\begin{proof}
Write $h(z)=\sum_{\substack{k=1\\P_kb=0}}^K\lambda_k\norm{P_kz}$ and let
$f$ be the gradient of the remaining penalty, including the ridge.
Put $r=q-f(b)$ and $C=\{v:h(v)=\langle r,v\rangle\}$.
Subtracting the differentiable summand from the supporting inequality
shows that every response blowup $g$ satisfies
\[
 g(z)\in C,\qquad
 \langle z-Hg(z),v-g(z)\rangle\le0\quad(v\in C).
\]
Lean proves this by passing the subgradient inequalities to the limit;
it does not assume a blowup formula.
Positive group weights give $C\cap(-C)=T$, and $C+T=C$.
The quadratic variational inequality yields
$g(z+Ht)=g(z)+t$ for $t\in T$ and
$\Pi_T Hg(z)=\Pi_Tz$.

If $g$ is linear on $F$, oddness implies $g(F)\subseteq T$.
The preceding translation identity extends this linearity to
$F+HT$. Maximal dimension forces $HT\subseteq F$, hence
$g(F)=T$. The projected equation identifies the curvature on this
range, and uniqueness of the symmetric factor gives the stated inverse.
\end{proof}

Lean checks the analytic limit in
\lean{groupLasso_response_blowup_quadratic}, maximality in
\lean{quadraticConeResponse_maximal_range}, and the resulting identities in
\lean{dataDrivenResponse_groupLasso_spec},
\lean{dataDrivenALOOperator_eq_groupLasso_reducedInverse},
\lean{dataDrivenALOOperator_eq_groupLassoALOOperator}, and
\lean{dataDrivenALOOperator_eq_groupLassoEmpiricalALOOperator}.
Zero-weight groups must be removed before using this recovery theorem:
the penalty function cannot detect such terms, but retaining them in
$T_{\rm grp}$ can artificially constrain the prescribed matrix.

\begin{theorem}[Replacement of the slice operator by the data-driven operator]
\label{thm:data-comparison}
Under the assumptions of \cref{thm:linear} or \cref{thm:single}, use the
single-index representation of the responses and define $A_{\rm si}$
as above. Write $s_i=\psi_{y_i}(x_i^T\widehat b)$ and $\alpha=n\mu$.
Almost surely on $\{s\ne0\}$,
\begin{equation}
 \rank(A_{\rm data}-A_{\rm si})\le4.
 \label{eq:data-rank}
\end{equation}
The following comparisons hold without an exceptional-event indicator:
\begin{align}
 \sum_{i=1}^n s_i^2\bigl\{\tr[(A_{\rm data}-A_{\rm si})\Sigma]\bigr\}^2
    &\le\frac{64}{n\mu^2}\quad\text{almost surely},
       \label{eq:data-trace-comparison}\\
 \E\sum_{i=1}^n s_i^2\bigl[W_i(A_{\rm data})-W_i(A_{\rm si})\bigr]^2
    &\le C_{\rm rank}(\mu,\delta).
       \label{eq:data-weight-comparison}
\end{align}
The random quantity in \eqref{eq:data-weight-comparison} is integrable.
When $s=0$, the weighted differences vanish, so no rank assertion
is needed on that event.
\end{theorem}

\begin{proof}
On each fixed signal slice, let $B,C$ be the derivatives of the
coefficient and its penalty subgradient in the original coordinates,
and let $\mathcal L_{\rm obs}=S^T\mathcal L$ be the transported forcing
map. The Lipschitz-score chain rule gives
$C=\mathcal L_{\rm obs}-MB$ and $A_{\rm si}(C+MB)=B$.
Every selected response blowup agrees with $B$ on directions $C\eta$.
If the score vector is nonzero, the codimension of $\ran C$ is at
most two. Maximality of the selected linearity subspace and symmetry
therefore give \eqref{eq:data-rank}.
Whitening the primal-dual relation preserves this rank bound.
It does not require the selection rule itself to commute with whitening.

Both whitened operators have norm at most $\alpha^{-1}$ almost
surely. Consequently the absolute trace of their difference is at
most $4(2/\alpha)$ on $\{s\ne0\}$. Since $|s_i|\le1$, this proves
\eqref{eq:data-trace-comparison}. The Gaussian finite-rank rational-weight
comparison proves \eqref{eq:data-weight-comparison}, using the global
norm and denominator bounds. Borel measurability permits averaging
the fixed-slice assertions and then integrating over the independent
noise. No fitted quantity is conditioned on in the conclusions.
\end{proof}

\paragraph{Lean formalization.}
The derivative transport and rank statements are
\lean{HasPenaltyModelJet.covariancePullback},
\lean{ae_covarianceSignalSlice_dataDriven_rank}, and
\lean{ae_covarianceDataDrivenSignalOperator_rank}.
Trace and weight control are provided by
\lean{sum_score_trace_difference_sq_le_of_finiteRank} and
\lean{signalOperator_weights_unconditional}.
The identity \lean{unwhiten_covarianceDataDrivenSignalOperator}
identifies the proof's whitened operator with $A_{\rm data}$ in the
observed coordinates. These are the comparisons used in the
fully instantiated endpoints of \cref{sec:unconditional}.

\begin{theorem}[Common derivative and low-rank discrepancy]
\label{thm:comparison}
Use the coordinates and forcing map of \cref{sec:operator}, and
write $\widetilde Q=SQS^T$ for a prescribed matrix $Q$.

\begin{enumerate}
\item \textbf{Twice continuously differentiable penalties, full design.}
Suppose $R$ is $C^2$ and the fixed losses $\ell_i$ are $C^2$, and take
$\mathcal E=\R^{n\times p}$.
At every differentiability point of $\widehat\theta$,
\begin{equation}
 B=\widetilde Q_{\rm sm}\mathcal L
    =\widetilde A\mathcal L.
 \label{eq:smooth-factor}
\end{equation}
If $v\ne0$, then
\begin{equation}
 \ran(\widetilde Q_{\rm sm}-\widetilde A)\subseteq\spn\{h\},
 \qquad \rank(Q_{\rm sm}-A)\le1.
 \label{eq:smooth-rank}
\end{equation}
In particular, exact equality $Q_{\rm sm}=A$ holds if, in addition,
\begin{equation}
 \text{there exists }i\text{ such that }v_i\ne d_i u_i.
 \label{eq:nonresonance}
\end{equation}
This condition, together with $v\ne0$, makes $\mathcal L$ surjective.
Without it, \eqref{eq:smooth-factor} asserts equality on the forcing
directions.

\item \textbf{Group Lasso, full design.}
For \eqref{eq:group-penalty} and fixed losses with 1-Lipschitz
derivatives, almost everywhere on the full design space,
\begin{equation}
 B=\widetilde Q_{\rm grp}\mathcal L
    =\widetilde A\mathcal L.
 \label{eq:group-factor}
\end{equation}
If $v\ne0$, then
\begin{equation}
 \ran(\widetilde Q_{\rm grp}-\widetilde A)\subseteq\spn\{h\},
 \qquad \rank(Q_{\rm grp}-A)\le1.
 \label{eq:group-rank}
\end{equation}
No strict inequality in the inactive-group KKT conditions and no
joint group coverage are required for these assertions.
If $\ran B=ST_{\rm grp}$, then $Q_{\rm grp}=A$ by uniqueness
of the symmetric factor.

\item \textbf{Group Lasso, single-index slice.}
Under the model and score conditions of \cref{thm:single} and the
group-ridge conditions in \eqref{eq:group-penalty}, put $s=S\beta$.
The two whitened matrices represent the same derivative along
$\mathcal E=\{E:Es=0\}$ almost surely.
On $\{v\ne0\}$ they satisfy
\begin{equation}
 \ran(\widetilde Q_{\rm grp}-\widetilde A_{\rm si})
       \subseteq\spn\{s,\widehat\theta\},\qquad
 \rank(Q_{\rm grp}-A_{\rm si})\le2.
 \label{eq:single-rank}
\end{equation}
The derivative and rank assertions do not require
\eqref{eq:coverage}, or a lower bound on the angle between
$s$ and $\widehat\theta$.
\end{enumerate}
The rank assertions are conditional on a nonzero fitted score
vector. If a row's score is zero, its ALO correction is zero for
either matrix, irrespective of their discrepancy.
\end{theorem}

\begin{proof}
For a $C^2$ penalty, differentiation of the KKT equation gives
$(\nabla^2\widetilde R+Z^TDZ)B=\mathcal L$, proving
\eqref{eq:smooth-factor}. For group Lasso, derivative locality in
inactive groups and the projected KKT equation yield
\eqref{eq:group-factor}; the same argument applies to valid
fixed-signal slices. The difference of the two symmetric matrices
annihilates the range of $\mathcal L$. When $v\ne0$, its range is
therefore contained in $\spn\{h\}$ on the full design space, or in
$\spn\{s,\widehat\theta\}$ on a single-index slice.
Surjectivity of $\mathcal L$ gives exact equality in the smooth
case. Equality of the derivative range and the prescribed active
space gives equality in the group case by
\eqref{eq:factor-definition}.
\end{proof}

\paragraph{Lean formalization.}
The smooth statements are
\lean{fderiv_eq_smoothEmpiricalALOOperator_forcing},
\lean{smoothEmpiricalALOOperator_sub_borel_range_le_span}, and
\lean{smoothEmpiricalALOOperator_eq_borel_of_nonresonance},
applied to the whitened penalty. Their current interfaces require
$C^1$ scores; this is why the smooth part of
\cref{thm:comparison} assumes $C^2$ losses even though the
statistical theorems require only Lipschitz scores.
The group statements are
\lean{ae_groupLasso_covariance_forcing},
\lean{ae_groupLasso_covariance_sub_borel_range_le_span},
\lean{ae_groupLassoSignal_difference_range_le}, and
\lean{ae_groupLassoSignal_difference_finrank_le_two}.
The exact-equality assertion for matching ranges combines
\lean{groupLassoALOOperator_spec} with
\lean{symmetric_factor_unique}.

\begin{corollary}[Unconditional bounds for the prescribed group inverse]
\label{cor:prescribed-group}
Use the group penalty \eqref{eq:group-penalty} and compute
$Q_{\rm grp}$ by \cref{def:group-inverse}.
Under the linear-regression assumptions of \cref{thm:linear},
$|\alo(Q_{\rm grp})-\loo|$ is integrable and
\begin{equation}
 \E|\alo(Q_{\rm grp})-\loo|
     \le \frac{C_{\rm grp,lin}(\mu,\delta)}{\sqrt n}.
 \label{eq:prescribed-linear}
\end{equation}
No joint group coverage is required here.
Under the single-index assumptions of \cref{thm:single}, if additionally
\begin{equation}
 \bigcap_{k=1}^K\ker P_k=\{0\},
 \label{eq:coverage}
\end{equation}
then the two quantities below are integrable and
\begin{align}
 \E\sum_{i=1}^n|t_i^{\rm ALO}(Q_{\rm grp})-t_i^{\rm LOO}|^2
    &\le C_{\rm grp,2}(\mu,\delta),\\
 \E|\alo(Q_{\rm grp})-\loo|
    &\le \frac{C_{\rm grp,si}(\mu,\delta)}{\sqrt n}.
\end{align}
For either model, the prescribed mean-field prediction also satisfies
\begin{align}
 \E\sum_{i=1}^n|t_i^{\rm MF}(Q_{\rm grp})-t_i^{\rm LOO}|^2
    &\le C_{\rm emp}(\mu,\delta),\\
 \E|\operatorname{MF}(Q_{\rm grp})-\loo|
    &\le \frac{C_{\rm grp,MF}(\mu,\delta)}{\sqrt n},
\end{align}
with both quantities integrable and without the coverage condition
\eqref{eq:coverage}. All constants depend only on $\mu,\delta$.
All expectations are under the full noise/design law.
For positive retained group weights, the identity in
\cref{thm:data-group-recovery} additionally gives all the bounds of
\cref{cor:data-group}, with its constants and without \eqref{eq:coverage}.
\end{corollary}

\paragraph{Lean formalization.}
The prescribed ALO endpoints are
\lean{groupLasso_prescribed_covarianceALO_test_bound_unconditional}
and \lean{singleIndex_groupLasso_ALO_unconditional}.
The prescribed mean-field endpoint is
\lean{singleIndex_groupLasso_meanField}; it includes regression by
the specialization $f_i(U,t)=t+\eps_i$.
The identity \lean{singleIndex_groupOperator_eq_original_loss}
identifies the prescribed matrix for the score-normalized losses with
that of the original losses.
The coverage assumption in the single-index prescribed ALO theorem
is a retained hypothesis of that older Lean endpoint; it is not required
by the data-driven group-penalty results in \cref{cor:data-group}, which
now apply to the prescribed matrix for positive retained weights.

\begin{remark}[Scope of the matrix identifications]
\Cref{thm:comparison} compares the derivative-based operator with
$Q_{\rm sm}$ or $Q_{\rm grp}$.
\Cref{thm:data-comparison} compares $A_{\rm data}$ with the
slice derivative-based operator. These are distinct comparisons.
Exact recovery $A_{\rm data}=Q_{\rm sm}$ and
$A_{\rm data}=Q_{\rm grp}$ (for positive retained group weights)
is established separately by
\cref{thm:data-smooth-recovery,thm:data-group-recovery}.
The proofs identify the response range and curvature; they do not
infer equality merely from a shared pseudoinverse form.

The general-penalty theorems in \cref{sec:unconditional} use
$A_{\rm data}$ throughout. \Cref{cor:prescribed-group} independently
uses $Q_{\rm grp}$ in both numerator and denominator. Its older
proof does not need the recovery identities. For $Q_{\rm sm}$, Lean also proves
\lean{sum_smoothEmpiricalALOPrediction_sub_borel_sq_le} on the row event
and the global bound \lean{smoothEmpiricalALOOperator_denominator};
the unconditional conclusions follow from the main theorems by the
pointwise identity in \cref{thm:data-smooth-recovery}.
\end{remark}

\section{Summary of the operator constructions and comparisons}
\label{sec:operator-summary}

Throughout this summary the matrices are in the original coefficient
coordinates, at the same empirical fit, with
$q=-X^Ts$, $M=X^TDX$, and
$G_R(u)=\argmin_b\{R(b)-u^Tb\}$.
The four constructions and the two prescribed matrices have the
following roles.

\begin{center}
\small
\begin{tabular}{@{}p{.17\linewidth}p{.40\linewidth}p{.35\linewidth}@{}}
\hline
Operator & Defining choice & Dependence and purpose \\
\hline
Canonical $A$ & Differentiate the fit on a prescribed design space;
factor its derivative through the forcing map
(\cref{def:operator}). & Depends on the perturbation space and offset.
The signal-slice instance $A_{\rm si}$ is an oracle proof operator. \\[4pt]
$A_{\rm data}$ & Select a subsequential penalty-response blowup and
a largest space on which it is linear
(\cref{def:data-response,def:data-operator}). & Uses only $R,q,M$;
ordinary differentiability is not required. \\[4pt]
$A_{\rm der}$ & Differentiate $G_R$ on a largest admissible dual
subspace (\cref{def:ordinary-penalty-response}). & Uses only $R,q,M$;
no subsequential limit is selected. \\[4pt]
$A_{\rm emp}$ & Differentiate the fixed-response empirical fit and
its subgradient on an admissible design space, maximizing forcing
rank (\cref{def:ordinary-design}). & Uses the observed losses, penalty
and design; the perturbation space is selected from these inputs. \\[4pt]
$Q_{\rm sm}$ & Invert $\nabla^2R(\widehat b)+M$
(\cref{def:group-inverse}). & Prescribed inverse for a $C^2$ penalty. \\[4pt]
$Q_{\rm grp}$ & Invert group curvature plus $M$ on the active-group
tangent space (\cref{def:group-inverse}). & Prescribed active-group
inverse, extended by zero. \\
\hline
\end{tabular}
\end{center}

The three replacements $A_{\rm data},A_{\rm der},A_{\rm emp}$ are
data-driven and require no covariance or signal parameter in their
construction. Their mean-field corrections still use $\Sigma$.
All have the structural form
\[
 A=\bigl[\Pi_T(\overline H+M)\Pi_T\bigr]^\dagger.
\]
This common formula does not identify the matrices: their selected
coefficient spaces $T$ and curvatures $H$ can differ.

\paragraph{Exact equalities.}
For a $C^2$ strongly convex penalty,
\begin{equation}
 A_{\rm data}=A_{\rm der}=Q_{\rm sm}.
 \label{eq:summary-smooth-equality}
\end{equation}
For the positive-weight group Lasso penalty with ridge,
\begin{equation}
 A_{\rm data}=A_{\rm der}=Q_{\rm grp}.
 \label{eq:summary-group-equality}
\end{equation}
These identities hold at every empirical fit, with no forcing-rank or
nonzero-score assumption. The group identity includes overlapping
groups and inactive boundary groups, without strict complementarity
or joint group coverage. See
\cref{thm:data-smooth-recovery,thm:data-group-recovery,sec:ordinary-penalty}.
More generally, ordinary differentiability of $G_R$ on the entire
dual space at $q$ makes the Section~\ref{sec:data-operator} and
Section~\ref{sec:ordinary-penalty} constructions agree, since every
response blowup then equals that full derivative.

If the full design space is admissible in
\cref{def:ordinary-design}, $A_{\rm emp}$ equals the full-design
canonical operator with the same fit, fixed losses and offset $c=0$.
Full-space admissibility alone does not imply
$A_{\rm emp}=A_{\rm data}$ or $A_{\rm emp}=A_{\rm der}$.
A sufficient additional condition is surjectivity of the empirical
forcing map on the selected space:
\begin{equation}
 \rank(\mathcal L_X|_{\mathcal E})=p
 \quad\Longrightarrow\quad
 A_{\rm emp}=A_{\rm der}=A_{\rm data}.
 \label{eq:summary-surjective-equality}
\end{equation}
Indeed, $\ker C=\ker(\mathcal L_X|_{\mathcal E})$, so $C$ is
surjective. The ordinary graph-jet argument in
\cref{thm:ordinary-penalty-statistics} makes $G_R$ differentiable on
the entire dual space. The three operators consequently solve the
same forcing equation, which determines the matrix when the forcing
map is surjective.

In general the constructions are distinct. The oscillatory example in
\cref{sec:ordinary-penalty} separates $A_{\rm data}$ from
$A_{\rm der}$. The zero-fit example in
\cref{sec:empirical-derivative} separates $A_{\rm emp}$ from both,
even for a smooth ridge penalty and full-space admissibility.

\paragraph{Low-rank comparisons.}
Under the hypotheses of \cref{thm:linear} or \cref{thm:single},
almost surely on $\{s\ne0\}$, the stated comparison theorems give
\begin{equation}
 \rank(A_\star-A_{\rm si})\le4,
 \qquad A_\star\in\{A_{\rm data},A_{\rm der},A_{\rm emp}\}.
 \label{eq:summary-slice-ranks}
\end{equation}
Here $A_{\rm si}$ is the particular signal-slice operator, not an
arbitrary choice of perturbation space in \cref{def:operator}.
For the prescribed matrices, \cref{thm:comparison} gives:
\begin{enumerate}
\item For a $C^2$ penalty and fixed $C^2$ losses, the full-design
canonical operator satisfies $\rank(Q_{\rm sm}-A)\le1$ on nonzero
scores, at differentiability points. The nonresonance condition
\eqref{eq:nonresonance} gives equality of the matrices.
\item With the group penalty, $\rank(Q_{\rm grp}-A)\le1$ for the
full-design canonical operator, almost everywhere on nonzero scores.
Equality follows if its coefficient derivative spans the
entire active-group tangent space.
\item On a single-index slice,
$\rank(Q_{\rm grp}-A_{\rm si})\le2$ almost surely on nonzero scores.
Thus \eqref{eq:summary-group-equality} improves the bound in
\eqref{eq:summary-slice-ranks} to two for both penalty-response
constructions in the group Lasso case.
\end{enumerate}

The bounds in \eqref{eq:summary-slice-ranks} immediately imply a
pairwise rank bound of eight. The common-response argument used in
their proofs gives the sharper consequence
\begin{equation}
 \boxed{\rank(A_1-A_2)\le4,
 \qquad A_1,A_2\in\{A_{\rm data},A_{\rm der},A_{\rm emp}\},}
 \label{eq:summary-pairwise-rank}
\end{equation}
almost surely on $\{s\ne0\}$. To see this, associate to each
construction its dual response space: the selected $F_q$ for
$A_{\rm data}$, the ordinary differentiability space $F$ for
$A_{\rm der}$, and $\ran C$ for $A_{\rm emp}$. Each has codimension
at most two, by the signal-slice competitor and the respective
maximality rule. On the intersection $W$ of any two such spaces,
the two linear responses agree, say with value $Kz$. This follows
from the ordinary graph expansion for the latter two constructions
and its agreement with every response blowup for the first.
Their factor equations imply
\[
 A_1(z+MKz)=Kz=A_2(z+MKz),\qquad z\in W.
\]
The map $z\mapsto z+MKz$ is injective on $W$: if it vanishes,
strong monotonicity gives
$a\norm{Kz}^2\le\langle Kz,z\rangle
=-\langle Kz,MKz\rangle\le0$, hence $Kz=z=0$.
Since $\operatorname{codim}W\le4$, rank--nullity proves
\eqref{eq:summary-pairwise-rank}. This pairwise consequence and
\eqref{eq:summary-surjective-equality} are recorded here with their
mathematical arguments; they are not separate assembled Lean endpoints.

The nonzero-score qualification is essential for the rank statements.
On zero-score samples, the matrices can differ in many directions,
but the score-weighted ALO and mean-field corrections vanish.
Full-space admissibility for $A_{\rm emp}$ is guaranteed almost surely
by \cref{prop:ordinary-full-space}, for example for countable-valued
responses, or for linear regression with independent noise having a
density. These conditions ensure admissibility, not forcing
surjectivity or equality of the matrices.

\paragraph{Statistical consequences.}
Under the common model, penalty and score hypotheses,
$A_{\rm data},A_{\rm der},A_{\rm emp}$ satisfy all four unconditional
ALO and mean-field conclusions with the same constants: expected
summed squared prediction errors are bounded, and expected absolute
test-criterion errors are of order $n^{-1/2}$. Their positivity,
covariance-relative norm and rational-denominator bounds permit the
finite-rank comparisons to be used quantitatively. Exact equality of
the matrices is not needed for these statistical conclusions.
By \eqref{eq:summary-smooth-equality} and
\eqref{eq:summary-group-equality}, the same conclusions apply to the
prescribed smooth and positive-weight group inverses. In particular,
the main group bounds require no coverage assumption; the coverage
condition in one older endpoint of \cref{cor:prescribed-group} is
unnecessary when using this exact-recovery route.

\end{document}